\documentclass[12pt]{iopart}

\usepackage{amsthm}
\usepackage{amssymb}
\usepackage{color}
\usepackage{graphics}
\usepackage{subfig}
\usepackage{ulem}
\usepackage{hyperref}
\usepackage{graphicx}
\usepackage{cleveref}
\usepackage{comment}

\newcommand{\R}{\mathbb{R}}

\newcommand{\N}{\mathbb{N}}

\newcommand*\Laplace{\mathop{}\!\mathbin\bigtriangleup}

\newcommand{\dist}{\mathrm{dist}}

\newcommand{\red}[1]{\textcolor{red}{#1}}

\newcommand{\fa}{\mathfrak{a}}
\newcommand{\fb}{\mathfrak{b}}
\newcommand{\fab}{\fa\fb}
\newcommand{\fba}{\fb\fa}
\newcommand{\muab}{\mu_{\fa\fb}}
\newcommand{\nuab}{\nu_{\fa\fb}}
\newcommand{\ha}{H^{\fa}}
\newcommand{\hb}{H^{\fb}}

\newcommand{\bv}{\mathbf{v}}

\newtheorem{thm}{Theorem}
\newtheorem{lem}{Lemma}
\newtheorem{cor}{Corollary}

\theoremstyle{definition}
\newtheorem{defn}{Definition}

\newtheorem{rem}{Remark}

\eqnobysec

\begin{document}
\title[Norm-dependent convergence and stability of the inverse scattering series]{Norm-dependent convergence and stability of the inverse scattering series for diffuse and scalar waves}

\author{Srinath Mahankali$^1$ and Yunan Yang$^2$}
\address{$^1$ Massachusetts Institute of Technology, Cambridge, MA 02139, USA}
 \address{$^2$ Institute for Theoretical Studies, ETH Z\"urich, Z\"urich, Switzerland 8092}
 \eads{\mailto{srinathm@mit.edu}, \mailto{yunan.yang@eth-its.ethz.ch}}
\vspace{10pt}
\begin{indented}
\item[]February 2023
\end{indented}

\begin{abstract}
This work analyzes the forward and inverse scattering series for scalar waves based on the Helmholtz equation and the diffuse waves from the time-independent diffusion equation, which are important PDEs in various applications. Different from previous works, which study the radius of convergence for the forward and inverse scattering series, the stability, and the approximation error of the series under the $L^p$ norms, we study these quantities under the Sobolev $H^s$ norm, which associates with a general class of $L^2$-based function spaces. The $H^s$ norm has a natural spectral bias based on its definition in the Fourier domain: the case $s<0$ biases towards the lower frequencies, while the case $s>0$ biases towards the higher frequencies. We compare the stability estimates using different $H^s$ norms for both the parameter and data domains and provide a theoretical justification for the frequency weighting techniques in practical inversion procedures. We also provide numerical inversion examples to demonstrate the differences in the inverse scattering radius of convergence under different metric spaces.
\end{abstract}

\section{Introduction}
We are interested in studying the forward scattering Born series and the inverse scattering Born series of both the Helmholtz equation and the diffuse wave equation. This work is motivated by the paper of Moskow and Schotland~\cite{moskow2008convergence}, where such studies were conducted under the $L^p$ function spaces for both parameter and measurement in which the radius of convergence and the Lipschitz stability constant of the inverse problem change with respect to the choice of $p$. There have been many works extending this analysis to various inverse problems, including the optical diffusion tomography~\cite{JS1}, diffuse waves~\cite{JS2}, scalar waves~\cite{JS5}, the electromagnetic scattering~\cite{JS3}, the Calder\'{o}n problem~\cite{JS4}, the Schr{\"o}dinger problem~\cite{bardsley2014restarted}, the radiative transport equation~\cite{JS6} and optimal tomography on  graphs~\cite{JS7}. In particular, studies of the radius of convergence and stability have been extended to general Banach spaces, and the parameter and data spaces do not have to be the same~\cite{moskow2019inverse}.
In this work, we set the parameter space and the measurement space as the $L^2$-based Sobolev space $\ha$ and $\hb$, respectively. We investigate how different choices of $\fa$ and $\fb$ could change the optimization behaviors of solving inverse problems based on the Helmholtz equation and diffuse wave equation, but similar strategies can apply to other inverse problems. The $L^2$-based Sobolev spaces are Hilbert spaces, and the inner products involve the Laplacian operator, thus making them easy to implement from a computational perspective. After discretization, the corresponding objective function becomes a weighted least-squares error, and the quadratic nature makes gradient computation efficient.

Practically speaking, inverse medium problems constrained by such equations are often formulated as PDE-constrained optimization problems, where well-posedness properties of the inverse problem are translated into the uniqueness of the minimizer, sensitivity to data noise, rate of convergence, and other features in an optimization framework. The choice of $\fb$ for the measurement space $\hb$ then corresponds to the choice of the objective function in the resulting nonlinear, nonconvex optimization problem. On the other hand, the choice of $\fa$ for the parameter space  $\ha$ imposes a priori information on the regularity of the parameter one aims to reconstruct. It could also lead to different gradient flow formulations when passing the Fr\'echet derivative to the gradient, thus giving rise to different gradient descent algorithms for solving such nonconvex optimization problems. Both choices affect the convergence rate and potentially change the stationary points to which the iterative gradient-based algorithm converges, even with the same initial guess~\cite{nurbekyan2022efficient}. We will demonstrate this later in~\sref{sec:numerical_inversion}.

There have been numerous works on changing the function spaces for both the parameter and the data in different applications. It is worth noting that the $H^{-1}$ semi-norm is closely related to the quadratic Wasserstein metric from optimal transportation~\cite{otto2000generalization,peyre2018comparison}. This connection has been utilized in many imaging and inverse problem applications~\cite{papadakis2014optimal,yang2018application,dunlop2020stability} and extended to the general $H^s$ norm as an objective function for data-matching inverse problems~\cite{engquist2020quadratic,yang2021implicit}. Methods based on the Sobolev gradient~\cite{neuberger2009sobolev}, where the gradient of a given functional with respect to the parameter function is taken with respect to the inner product induced by the underlying Sobolev norm, have demonstrated advantages in image sharpening and edge-preserving~\cite{calderMY10,sundaramoorthi2007sobolev}. In~\cite{jacobs2019solving}, it was shown that choosing $H^1$ gradient flow for minimizing $L^1$-type objective functions yields a Lipschitz constant independent of the discretization grid size, achieving an optimal rate of convergence in the gradient descent method. In~\cite{ren2022generalized}, the impacts of changing both the data and parameter spaces were discussed in machine learning in a discrete setup.

To tackle the local minima issues due to the inherent nonconvexity and sensitivity to noise of the $L^2$ norm as the objective function, 
there is a general strategy of adding adaptive weights to different frequency components of the data, which is often referred to as \textit{Frequency Marching}~\cite{Bao2015,Barnett2017},  \textit{Multiscale/Hierarchical Inversion}~\cite{bunks1995multiscale,Fichtner2013} or \textit{Frequency Weighting} techniques for nonlinear inverse problems~\cite{Brossier2010,Oh2013,Kamei2014, jeong2017comparison}. The model parameters should be updated along directions such that both low-wavenumber and high-wavenumber structures can be appropriately resolved. One way to control the order of updating the low- and high-wavenumber model components is through the choice of $H^\fa$ parameter space. Since they are also influenced by the low- and high-frequency components of the data residual through the gradient descent update,  it is also preferable to have adaptive weights regarding different data frequencies in the optimization, which can be achieved by tuning different $\fb$ for the data metric space $H^{\fb}$ in our framework.

The main novelty of this work is to utilize the natural frequency biasing features of the $L^2$-based Sobolev spaces, as briefly mentioned above, and prove the impact of the choice of the parameter and data function spaces on the properties of the inverse scattering problems. We also focus on the class of $L^2$-based Sobolev space for its particular convenience in implementation, while we believe similar conclusions could be drawn for other classes of Sobolev spaces $W^{k,p}$ with $p$ fixed and $k\in\mathbb{Z}$. To the best of our knowledge, utilizing the spectral properties of the $L^2$-based Sobolev space to improve the radius of convergence and stability of the inverse Born series has yet to be considered in the literature. There are three main contributions in our work. First, we rigorously study the convergence property of the forward scattering series that maps the parameter in $H^\fa$ to the scattering data in $H^\fb$, and the inverse scattering series that maps from $H^{\fb}$ to $H^{\fa}$. We show that the radius of convergence for the forward scattering series is only affected by $\fb$, but the radius of convergence for the inverse scattering series is affected by both $\fa$ and $\fb$. Second, through a sequence of theorems, we demonstrate that the choice of $(\fa,\fb)$ consequently changes the stability of the limit of the inverse scattering series with respect to the small changes in the observed scattering data and changes the approximation error of the inverse scattering series. Our results show that one can improve the stability by tuning $(\fa,\fb)$. Third, we demonstrate the impact of the choice of these two function spaces in a few numerical inversion examples to illustrate how one can qualitatively change the performance of an inversion algorithm.

The rest of the paper is organized as follows. We first present some essential background in~\sref{sec:background}, where we briefly review the problem setups for the Helmholtz equation and the diffuse wave equation and present assisting lemmas that will be used in the later analysis. In~\sref{sec:forward_series}, we analyze the radius of convergence of the forward scattering series while embedding the PDE solution into $H^\fb$  and the variable coefficient into $H^\fa$. 
It is followed by~\sref{section:inverse_scattering_series} where we study the inverse scattering operator that maps $H^\fb \rightarrow H^{\fa}$ regarding its convergence, stability, and approximation error. There has been a recent result~\cite{hoskins2022analysis} deriving the inverse scattering radius of convergence using the geometric function theory under much weaker assumptions. We also incorporate analogous results here for our $H^s$-type function spaces. In~\sref{sec:discussions}, we discuss the results from sections~\ref{sec:forward_series} and~\ref{section:inverse_scattering_series} and how the important quantities change under different choices of $(\fa,\fb)$ pairs. We also compare the radius of convergence and approximation error obtained from the classic approach in~\cite{moskow2008convergence} and the geometric approach~\cite{hoskins2022analysis}. It is followed by numerical illustrations of the physical differences between the scalar wave and the diffuse wave equations on the radius of convergence for both series. In~\sref{sec:numerical_inversion}, we show some numerical inversion examples to demonstrate the improved stability and radius of convergence under proper $(\fa,\fb)$ choices. Conclusion follows in~\sref{sec:conclusion}.

\section{Preliminaries}\label{sec:background}
In this section, we review some background on the forward problems and some preliminary results for the $L^2$-based Sobolev norms. 

\subsection{Diffuse and Scalar Waves}
We consider the propagation of a scalar wave into a medium whose pressure field $u$ satisfies the Helmholtz equation
\begin{equation}\label{eq:helm-1}
    \Laplace_x u + \frac {\omega^2}{c(x)^2}u(x) = -S(x-x_1,\omega),\quad x\in\R^n,
\end{equation}
where $\omega$ is the frequency, $S(x,\omega)$ is the source term, $x_1$ is the position of the source term, and $c(x)$ is the spatial-dependent velocity. \Eref{eq:helm-1} can be derived from the scalar wave equation assuming the solution is time harmonic. It governs the propagation of time harmonic acoustic waves of small amplitude in a slowly varying inhomogenous medium~\cite{colton1998inverse}. 

We assume the inhomogenous region is contained inside a ball $B_a$ of radius $a$, i.e., $c(x) = c_0$ for $x\in\R^n \setminus B_a$. Let the wavenumber $k = \omega/c_0$ and $\eta(x) = c_0^2 / c^2(x) - 1$. The scattering problem becomes 
\begin{equation} \label{eq:helm-2}
    \Laplace_x u + k^2 (1+\eta) u = -
    S(x-x_1,\omega),\quad x\in \R^n,
\end{equation}
where $u = u_i + u_s$ with the incident wavefield $u_i$ satisfying 
\begin{equation*}
    \Laplace_x u_i + k^2 u_i =   -S(x-x_1,\omega),\quad x\in \R^n.
\end{equation*}
Here, the scattering wavefield $u_s$ satisfies 
\[
  \Laplace_x u_s + k^2 u_s =   -k^2 \eta(x) u(x,x_1),
\]
and the Sommerfeld radiation condition
\begin{equation} \label{eq:sommerfeld}
    \lim_{r\to \infty} r^{\frac{n-1}{2}}\left( \frac{\partial u_s}{\partial r} - iku_s\right) = 0, \qquad u_s(x,x_1) = u_s(r\hat{x},x_1),
\end{equation}
where $\|\hat{x}\| = 1$. Note that we assume $\eta(x)\geq -1$ for all $x\in B_a$ and $\eta(x)= 0$ outside $B_a$. The solution to \eref{eq:helm-2} satisfies the Lippmann--Schwinger equation
\begin{equation}\label{eq:LS}
    u(x,x_1) = u_i(x) + k^2 \int_{\R^n} G(x,y,k) u(y,x_1) \eta(y) dy ,
\end{equation}
where $u_i(x) = \int_{\R^n} G(x,y,k) S(y,\omega) dy$ and $G(x,y,k)$ is the Green's function such that with $y$ fixed,
\[
        \Laplace_x G + k^2 G = -\delta(x-y).
\]
Note that $G(x,y,k) = G(y,x,k)$. The expression of the Green's function depends on the dimensionality of the problem. For example, based on the outgoing Sommerfeld radiation condition~\eref{eq:sommerfeld},
\begin{equation}\label{eq:greens_helmholtz}
    G(x,y,k) = \cases{
    \frac{i\exp(ik|x-y|)}{2k},& $n = 1$, \\
    \frac{i}{4} H_0^{(1)} (k|x-y|),& $n = 2$, \\
    \frac{\exp(ik|x-y|)}{4\pi |x-y|},& $n = 3$. 
    }
\end{equation}
Here, $H_0^{(1)} $ denotes the Hankel function  of the first kind
and $|x-y|$ is the Euclidean distance between $x,y \in \R^n$.  

We can also apply our analysis to the diffuse wave equation~\cite{mandelis2000diffusion,moskow2008convergence},
\begin{equation}\label{eq:diffuse eqn}
 \Laplace_x u -k^2(1+\eta) u(x) = -S(x-x_1,\omega),
\end{equation}
where $u$ represents the energy density and $\eta(x)$ is the absorption coefficient, again assumed to be compactly supported on $B_a$. We refer to~\cite[Sec.~2]{moskow2008convergence} for details regarding the forward problem. Let $G(x,y,k)$ be the Green's function where
\begin{equation*}
    \Laplace_x G - k^2 G = -\delta(x-y),
\end{equation*}
for a fixed $y$ and $k$. The analytic form of $G$ also depends on the dimension~\cite{mandelis2013diffusion},
\begin{equation}\label{eq:greens_diffuse}
    G(x,y,k) = \cases{
    \frac{\exp(-k|x-y|)}{2k}, & $n = 1$, \\
    \frac{1}{2\pi}K_0(k|x-y|), & $n = 2$, \\
    \frac{\exp(-k|x-y|)}{4\pi |x-y|}, & $n = 3$,
    }
\end{equation}
where $K_0$ has the following integral representations for $x > 0$
\begin{equation}\label{eq:2d_diffuse_greens_function}
    K_0(x) = \int_{1}^{\infty} \frac{\exp(-xt)}{(t^2 - 1)^{\frac{1}{2}}} dt = \int_{0}^{\infty} \exp(-x\cosh t) dt.
\end{equation}

We can also define the incident and scattering wavefields for the diffuse wave equation~\eref{eq:diffuse eqn} under similar assumptions to the case for the Helmholtz equation~\eref{eq:helm-2}. Later, we will use the same set of notations for both the diffuse and scalar waves, except they have different Green's functions.

\subsection{The $L^2$-Based Sobolev Space}
The following definition of the Sobolev norm is based on the Sobolev space $W^{s,p}(\Omega)$ for nonnegative $s$~\cite{arbogast1999methods}.
\begin{defn}[Sobolev Space $W^{s,p}(\Omega)$]
Let $1 \le p < \infty$ and $s$ be a nonnegative integer. If a function $f$ and its weak derivatives $D^\alpha f = \frac{\partial^{|\alpha|}f}{\partial x_1 ^{\alpha_1}\cdots \partial x_n^{\alpha_n}}, |\alpha| \le s$ all lie in $L^p(\Omega)$, where $\alpha$ is a multi-index and $|\alpha| = \sum_{i=1}^n \alpha_i$, we say $f \in W^{s,p}(\Omega)$ and define the $W^{s,p}(\Omega)$ norm of $f$ as 
\begin{equation*}
\|f\|_{W^{s,p}(\Omega)} := \bigg(\sum_{|\alpha| \le s} \|D^{\alpha}f\|_{L^p(\Omega)}^p\bigg)^{\frac{1}{p}}.
\end{equation*}
\end{defn}
We also define the space $W^{s,p}_0(\Omega)$ as the space of functions $f \in W^{s,p}(\Omega)$ with compact support~\cite{arbogast1999methods}. Next, we define the Sobolev norm on the boundary $\partial\Omega$, where the set $\Omega$ is bounded with sufficiently smooth boundaries. We use the following definition from~\cite[p.4]{necas2011direct}.
\begin{defn}[Bounded Domains with Continuous Boundary]\label{defn:continuous_boundary}
Let $\Omega \subset \R^n$ be a bounded domain with the boundary $\partial\Omega$. Suppose the following conditions are satisfied:
\begin{itemize}
    \item There exist $\lambda,\gamma > 0$, and systems of local charts $(x_{r1},x_{r2},\ldots,x_{rn}) = (x_r',x_{rn})$ for $r \in \{1,2,\ldots,m\}$, and continuous functions $a_r$ defined on the closed $n-1$ dimensional cubes defined by $|x_{ri}| \le \lambda$ for $i \in \{1,2,\ldots,n-1\}$, such that every point $x \in \partial \Omega$ can be represented as $(x_r',a(x_r'))$ for at least one $r \in \{1,2,\ldots,m\}$.
    \item Let $\overline{\Delta_r}$ denote the set of points $x_r'$ such that $|x_{ri}| \le \lambda$ for $i \in \{1,2,\ldots,n-1\}$. The points $(x_r',x_{rn})$ where $x_r' \in \overline{\Delta_r}$ and $a_r(x_r') < x_{rn} < a_r(x_r') + \gamma$ are in $\Omega$, while the the points $(x_r',x_{rn})$ where $x_r' \in \overline{\Delta_r}$ and $a_r(x_r') - \gamma < x_{rn} < a_r(x_r')$ are not in $\overline{\Omega}$. 
\end{itemize}
Then, the boundary $\partial\Omega$ is called \textit{continuous}. Furthermore, if the functions $a_r$ are all Lipschitz, then we say $\Omega$ has Lipschitz boundary $\partial\Omega.$
\end{defn}

For example, the systems of local charts for the surface of a unit cube in $\mathbb{R}^3$ are defined on the six faces. Next, we use the definition for $\mathfrak{R}^{k,\mu}$ domains from~\cite[p.49]{necas2011direct}.

\begin{defn}[Bounded Domains of Type $\mathfrak{R}^{k,\mu}$]
Let $k$ be a nonnegative integer or infinity, and let $0 \le \mu \le 1$. Consider a bounded domain $\Omega$ and suppose there exist functions $a_r,$ with $r \in \{1,2\ldots,m\}$ defined on $\overline{\Delta_r},$ with continuous derivatives of order $\le k$.
Furthermore, suppose for all multi-indices $\alpha$ with $|\alpha| \le k$ and for any $x_r', y_r' \in \overline{\Delta_r}$, $|D^\alpha a_r(x_r') - D^\alpha a_r(y_r')| \le c |x_r' - y_r'|^\mu$, for some constant $c > 0$. Then, we say $\Omega \in \mathfrak{R}^{k,\mu}$.
\end{defn}

Using the above definitions, we can define Sobolev spaces on $\partial\Omega$; see~\cite[p.83]{necas2011direct}.
\begin{defn}[Sobolev Space $W^{s,2}(\partial\Omega)$]\label{defn:sobolev_w_space_boundary}
Let $s$ be a positive integer and consider a domain $\Omega \in \mathfrak{R}^{s-1,1}$. We define the space $W^{s,2}(\partial\Omega)$ as the space of functions $f$ for which $f(x_r',a_r(x_r')) = f_r \in W^{s,2}(\Delta_r)$. We define the $W^{s,2}(\partial\Omega)$ norm of $f$ as
\begin{equation*}
    \|f\|_{W^{s,2}(\partial\Omega)} := \bigg(\sum_{r=1}^{m} \|f_r\|_{W^{s,2}(\Delta_r)}^2\bigg)^{\frac{1}{2}}.
\end{equation*}
The space $W^{s,2}(\partial\Omega)$ also has an inner product: for functions $f,g \in W^{s,2}(\partial\Omega)$, we define 
\begin{equation}\label{eq:boundary H^s}
    \langle f,g \rangle_{W^{s,2}(\partial\Omega)} := \sum_{r=1}^{m} \langle f_r, g_r \rangle_{W^{s,2}(\Delta_r)}
\end{equation}
where for each $r \in \{1,2,\ldots,m\}$, $\langle f_r, g_r \rangle_{W^{s,2}(\Delta_r)} = \sum_{|\alpha| \le s}\int_{\Delta_r} D^\alpha f_r(x_r')\, \overline{D^\alpha g_r(x_r')} dx_r'$.
\end{defn}

We also extend~\Cref{defn:sobolev_w_space_boundary} to the case where $s$ is a negative integer, by defining the space $W^{-s,2}(\partial\Omega)$ as the dual of the space $W^{s,2}(\partial\Omega)$, for positive integers $s$.
\begin{defn}[Sobolev Norm $\|\cdot\|_{W^{-s,2}(\partial\Omega)}$]\label{defn:negative_sobolev_norm_boundary}
Consider a domain $\Omega$ with continuous boundary $\partial\Omega$. For functions $f$ and $g$ on $\partial\Omega$, let $\langle f,g \rangle =    \langle f,g \rangle_{W^{0,2}(\partial\Omega)}$ following~\eref{eq:boundary H^s}.
For a function $f$, we define its $W^{-s,2}(\partial\Omega)$ norm as
\[
    \|f\|_{W^{-s,2}(\partial\Omega)} := \sup_{g \in W^{s,2}(\partial\Omega)} \frac{\langle f,g \rangle}{\|g\|_{W^{s,2}(\partial\Omega)}}.
\]
If $\|f\|_{W^{-s,2}(\partial\Omega)} < \infty$, we say $f \in W^{-s,2}(\partial\Omega)$.
\end{defn}
For convenience, we will refer to the space $W^{s,2}$ as $H^s$ and $W^{s,2}_0$ as $H^s_0$ through the rest of this paper. Next, we state the following Poincar\'e's inequality from~\cite[Eqn.~(18.1)]{rektorys2012variational}, with a proof similar to~\cite[Theorem 12.17]{leoni2017first}.
\begin{thm}[Poincar\'e's Inequality in $H^1_0(B_a)$]\label{thm:poincare_inequality}
For any $f \in H^1_0(B_a)$, where $B_a \subset \R^n$ is a ball of radius $a$, and for any $1 \le i \le n$, we have
\[
    \int_{B_a}|f(x)|^2dx \le 2a^2\int_{B_a}\left|\frac{\partial}{\partial x_i}f(x)\right|^2 dx.
\]
\end{thm}

If $f \in H^s_0(B_a)$ for a positive integer $s\geq 1$, \Cref{thm:poincare_inequality} allows us to obtain a lower bound on $\|f\|_{H^s(B_a)}$ in terms of $\|f\|_{L^2(B_a)}$, which leads to the following~\Cref{lem:applied_poincare_inequality}.
\begin{lem}\label{lem:applied_poincare_inequality}
Let $f \in H^s_0(B_a)$, where $s$ is a positive integer and $B_a \subset \R^n$. We define
\begin{equation}\label{eq:psan}
    P(s,a,n) := \left(\sum_{j = 0}^{s} \frac{{n+j-1 \choose n-1}}{(2a^2)^j}\right)^{-\frac{1}{2}}.
\end{equation}
Then, $\|f\|_{L^2(B_a)} \le P(s,a,n)\|f\|_{H^s(B_a)}$.
\end{lem}
\begin{proof}
Recall that
\[
    \|f\|_{H^s(B_a)}^2 = \sum_{|\alpha| \le s} \|D^{\alpha}f\|_{L^2(B_a)}^2 =\sum_{j = 0}^{s}\sum_{|\alpha| = j} \|D^{\alpha}f\|_{L^2(B_a)}^2.
\]
For a particular  multi-index $\alpha$ where $|\alpha| = j$, after applying~\Cref{thm:poincare_inequality} for a total of $j$ times, we have 
\[
    \frac{1}{(2a^2)^j}\int_{B_a}|f(x)|^2dx \le \int_{B_a}\left|D^{\alpha}f(x)\right|^2 dx.
\]
This gives 
\begin{eqnarray*}
    \|f\|_{H^s(B_a)}^2  &\ge \sum_{j = 0}^{s}\sum_{|\alpha| = j} \frac{1}{(2a^2)^j}\|f\|_{L^2(B_a)}^2 = \|f\|_{L^2(B_a)}^2 \sum_{j = 0}^{s}\sum_{|\alpha| = j} \frac{1}{(2a^2)^j}.
\end{eqnarray*}
The number of multi-indices $\alpha$ s.t.~$|\alpha| = j$ is the number of $n$-tuples $(\alpha_1,\ldots,\alpha_n)$ of nonnegative integers ($\alpha_i\in\mathbb{N}$, $\forall i$) s.t.~$\alpha_1 + \cdots + \alpha_n = j$. Thus, there are ${n+j-1 \choose n-1}$ multi-indices $\alpha$ with $|\alpha| = j$, i.e., $\sum_{|\alpha| = j}  1 = {n+j-1 \choose n-1}$, which completes the proof.
\end{proof}
We will prove the following lemma as another generalization of~\Cref{thm:poincare_inequality}.
\begin{lem}\label{lem:generalized_poincare}
Let $f \in H^{s+1}_0(B_a)$, where $s$ is a positive integer and $B_a \subset \R^n$. Then, the following inequality holds:
\begin{equation*}
    \frac{1}{{n+s-1 \choose n-1}} \sum_{|\alpha| = s}\|D^{\alpha}f\|_{L^2(B_a)}^2 \le \frac{2a^2}{{n+s \choose n-1}} \sum_{|\beta| = s+1}\|D^{\beta}f\|_{L^2(B_a)}^2 .
\end{equation*}
\end{lem}
\begin{proof}
Consider pairs of multi-indices $\alpha, \beta$ where $|\alpha| = s$, $|\beta| = s+1$, and $\alpha \leq \beta$. Thus, $\exists i\in \{1,\ldots,n\}$ where $\beta_i = \alpha_i +1$. We define a ``weight" function $w(\alpha,\beta) := \beta_i / |\beta|$ on such pairs of multi-indices. For a fixed $\beta$ where $|\beta|=s+1$, we have
\begin{equation}\label{eq:w_split_alpha}
    \sum_{|\alpha|=s,\alpha\leq \beta} w(\alpha,\beta) = \sum_{\{i : \beta_i > 0\}} \frac{\beta_i}{s+1} = 1,
\end{equation}
where the sum is taken over all valid $\alpha$. Furthermore, for a fixed $\alpha$ where $|\alpha|=s$,
\begin{equation}\label{eq:w_split_beta}
    \sum_{{|\beta|=s+1,\alpha\leq \beta}} w(\alpha,\beta) 
    = \sum_{i = 1}^{n} \frac{\alpha_i+1}{s+1} = \frac{n+s}{s+1}.
\end{equation}
Using \Cref{thm:poincare_inequality} and \eref{eq:w_split_alpha}-\eref{eq:w_split_beta}, we have the following for $\alpha\leq \beta$,
$$
   2a^2 \sum_{|\beta| = s+1}\|D^{\beta}f\|_{L^2(B_a)}^2 \ge  \sum_{|\beta| = s+1} \sum_{|\alpha|=s} w(\alpha,\beta) \|D^{\alpha}f\|_{L^2(B_a)}^2
    =  \sum_{|\alpha| = s} \frac{n+s}{s+1} \cdot \|D^{\alpha}f\|_{L^2(B_a)}^2.
$$
This gives the desired inequality after dividing both sides by ${n+s \choose n-1}$.
\end{proof}

Next, we prove the following~\Cref{lem:greens_function_c_infty}, which will allow us to apply~Definitions~\ref{defn:sobolev_w_space_boundary} and~\ref{defn:negative_sobolev_norm_boundary} to the Green's function of both the Helmholtz and diffuse wave equations, under the condition that $\dist(\partial\Omega, B_a) = \sup_{x\in \partial\Omega,\, y\in B_a} |x-y| > 0$.

\begin{lem}\label{lem:greens_function_c_infty}
Let $s\in\mathbb{N}^+$ and the domain $\Omega \in \mathfrak{R}^{s-1,1}$.
Suppose $\dist(\partial \Omega, B_a) = \epsilon > 0$. For any $C^\infty$ function $f: [\epsilon,\infty) \mapsto \mathbb{R}$, $\widetilde{f}(x)  := f(|x-y|) \in H^s(\partial\Omega)$ for any fixed $y \in B_a$.
\end{lem}

\begin{proof}

We can express $\widetilde{f}$ as the composition of two functions $f \circ g$, where $g(x) = |x-y|$. First, we claim that $g \in H^s(\partial\Omega)$. For a fixed local chart $\Delta_r$, we can write $g(x)$ as
\[
    g_r(x_r') = |(x_r',a_r(x_r')) - y| = \sqrt{|x_r' - y'|^2 + (a_r(x_r') - y_n)^2},
\]
where $y' = (y_1,\ldots,y_{n-1})$. Since $\Omega \in \mathfrak{R}^{s-1,1}$, by definition $a_r \in C^{s-1}(\Delta_r)$. Furthermore, the derivatives of $a_r$ of order at most $s-1$ are all Lipschitz. Thus, from Rademacher's theorem~\cite{federer2014geometric}, the derivatives of $a_r$ of order $s$ exist almost everywhere and we get that $a_r \in H^s(\Delta_r)$. We also have $g_r(x_r') \in H^s(\Delta_r)$ based on the chain rule. 

We define $\widetilde{f}_r(x_r') = \widetilde{f}(x_r',a_r(x_r')) = f(g_r(x_r'))$
on the local chart $\Delta_r$. Since the range of $g_r(x_r')$ is a subset of $[\epsilon,\infty)$, and $f \in C^\infty([\epsilon,\infty))$, we have $\widetilde{f}_r(x_r') \in H^s(\Delta_r)$ by the chain rule. Since for every local chart $\Delta_r$, we have $\widetilde{f}_r(x_r') \in H^s(\Delta_r)$, we get that $\widetilde{f}(x) \in H^s(\partial\Omega)$.
\end{proof}

\begin{rem}
We can apply \Cref{defn:sobolev_w_space_boundary} to the Green's function of the 2D and 3D Helmholtz equations (see~\eref{eq:greens_helmholtz}), since the necessary conditions to use \Cref{lem:greens_function_c_infty} apply to such $G(x,y)$. In 2D, the Green's function for the Helmholtz equation is $G(x,y) = \frac{i}{4} H_0^{(1)} (k|x-y|) = \frac{i}{4} H_0^{(1)} (kr)$, where $H_0^{(1)}(x) \in C^\infty((0,\infty))$ is the so-called Hankel function, and thus, $\frac{i}{4} H_0^{(1)} (kr) \in C^{\infty}([\epsilon,\infty))$ for any fixed $\epsilon > 0$. In 3D, the Green's function is $G(x,y) = f(|x-y|)$, where $f(r) = \exp(ikr)/(4\pi r)$. From the product rule, $f^{(n)}(r)$ is a polynomial in $\exp(ikr)$ and $r^{-1}$, with complex coefficients, $\forall n\in \N$. Thus, $\forall \epsilon > 0$, we also have that $f(r) \in C^{\infty}([\epsilon,\infty))$.

Similarly, we can also apply \Cref{defn:sobolev_w_space_boundary} to the Green's function of the 2D and 3D diffuse wave equations (see~\eref{eq:greens_diffuse}).  In 2D, the Green's function is  $G(x,y) = \frac{1}{2\pi}K_0(k|x-y|) = \frac{1}{2\pi}K_0(kr)$, where $K_0(x)$ has the integral representation in \eref{eq:2d_diffuse_greens_function} for $x > 0$.
From \eref{eq:2d_diffuse_greens_function}, we get $K_0^{(n)}(x) = \int_{0}^{\infty} (-\cosh t)^n \exp(-x\cosh t) dt$. Thus, $K_0(x) \in C^{\infty}((0,\infty))$. This implies that for any fixed $\epsilon > 0$, we must have $\frac{1}{2\pi}K_0(kr) \in C^{\infty}([\epsilon,\infty))$. Finally, in 3D, the Green's function is $G(x,y) = f(|x-y|)$, where $f(r) = \exp(-kr)/(4\pi r)$. From the product rule, $f(r) \in C^{\infty}([\epsilon,\infty))$ for any fixed $\epsilon > 0$.
\end{rem}

We present~\Cref{lem:linear_map_sobolev_norm_boundary}, which will be used in \sref{sec:forward_series} to obtain bounds for $\|f\|_{H^s(\partial\Omega)}$.
\begin{lem}\label{lem:linear_map_sobolev_norm_boundary}
Given $s\in \mathbb{Z}$, let $d$ be the total number of different multi-indices $\alpha$ where $|\alpha| \le |s|$. There exists a linear operator
\[
    \mathbf{T} : H^s(\partial\Omega) \to L^2(\partial\Omega)\times\cdots\times L^2(\partial\Omega)
\]
such that $\mathbf{T}f:\partial \Omega \mapsto \R^d$ and $\|f\|_{H^s(\partial\Omega)} = \|\mathbf{T}f\|_{L^2}$
for all functions $f \in H^s(\partial\Omega)$.
\end{lem}

\begin{proof}
First, if $s = 0$, then we can take $\mathbf{T}$ to be the identity. Next, consider the case where $s > 0$.  We define the linear operator $\mathbf{T}$ such that for a point 
\[
    x = (x_{r1},x_{r2},\ldots,x_{rn}) = (x_r',x_{rn}) = (x_r',a_r(x_r'))
\]
as defined in \Cref{defn:continuous_boundary}, and a function $f \in H^s(\partial\Omega)$, $\mathbf{T}f(x)$ is the vector with components $D^\alpha f_r(x_r')$ for each multi-index $\alpha \in \R^{n-1}$ with $|\alpha| \le s$. With this definition of $\mathbf{T}$, we note that $\mathbf{T}f$ is a piecewise function which depends on the local chart containing the input $x = (x_r',a_r(x_r'))$. However, $\mathbf{T}$ is still a linear operator in $f$. From \Cref{defn:sobolev_w_space_boundary}, this gives $\|f\|_{H^s(\partial\Omega)} = \|\mathbf{T}f\|_{L^2}$. 

Also, if $s > 0$, $\|\cdot\|_{H^{-s}(\partial\Omega)}$ is defined following~\Cref{defn:negative_sobolev_norm_boundary}. For fixed $f \in H^{-s}(\partial\Omega)$, and any $g \in H^{s}(\partial\Omega)$, we have
$\langle f,g \rangle \le\|f\|_{H^{-s}(\partial\Omega)}\|g\|_{H^s(\partial\Omega)}$.
Hence, $\langle f,g \rangle$ is a bounded linear functional on $H^s(\partial\Omega)$ for a fixed $f$. From the Riesz representation theorem, there is a unique $u \in H^s(\partial\Omega)$ such that for all $g \in H^s(\partial\Omega)$, $\langle f,g \rangle = \langle u,g \rangle_{H^s(\partial\Omega)}$,
and moreover, $\|f\|_{H^{-s}(\partial\Omega)} = \|u\|_{H^s(\partial\Omega)}$. Since $u$ is unique, there is an operator $\mathcal{L}: H^{-s}(\partial\Omega) \to H^s(\partial\Omega)$ such that $u = \mathcal{L}f$. 
It is easy to check that $\mathcal{L}$ is a linear operator and $\|f\|_{H^{-s}(\partial\Omega)}  = \|\mathcal{L}f\|_{H^s(\partial\Omega)}$. From the proof for the case $s \ge 0$, there exists a linear operator $\mathbf{T}_1$ such that $\|u\|_{H^s(\partial\Omega)} = \|\mathbf{T}_1 u\|_{L^2(\partial\Omega)}.$ Thus,
\[
    \|f\|_{H^{-s}(\partial\Omega)} = \|u\|_{H^s(\partial\Omega)} = \|\mathbf{T}_1 u\|_{L^2(\partial\Omega)} = \|\mathbf{T}_1 \mathcal{L} f\|_{L^2(\partial\Omega)},
\]
giving $\mathbf{T} = \mathbf{T}_1 \mathcal{L}$ as the desired linear operator.
\end{proof}

\section{Forward Scattering Series}\label{sec:forward_series}

Similar to~\cite[Eqn.\ (7)]{moskow2008convergence}, we apply fixed point iteration beginning with $u_i$ in~\eref{eq:LS}, which gives the following infinite series:
\begin{eqnarray}
    u(x,x_1) &=& u_i(x,x_1) + k^2 \int_{B_a} G(x,y_1) \eta(y_1) u_i(y_1,x_1) dy_1 + \label{eq:forward_series} \\
           &&  k^4 \int_{B_a}\int_{B_a} G(x,y_1)\eta(y_1) G(y_1,y_2)\eta(y_2) u_i(y_2,x_1) dy_1dy_2 + \ldots.\nonumber
\end{eqnarray}
Letting $\phi = u - u_i$, we can express the infinite series in the following manner:
\begin{equation}\label{eq:original_forward_series}
    \phi = K_1\eta + K_2\eta\otimes\eta + K_3\eta\otimes\eta\otimes\eta + \cdots
\end{equation}
where
\begin{equation}\label{eq:K_operator_formula} \fl
    (K_j f)(x,x_1) = k^{2j}\int_{\mathcal{B}_j}G(x,y_1)G(y_1,y_2)\cdots G(y_{j-1},y_j) u_i(y_j,x_1)\, f\,dy_1\cdots dy_j,
\end{equation}
where for $j\ge 1$,
\begin{eqnarray}
    \mathcal{B}_j &=& B_a\times \cdots\times B_a, \label{eq:B def}\\
    f &=& f(y_1,\ldots,y_j) = \eta(y_1)\cdots\eta(y_j). \label{eq:f def}
\end{eqnarray}

\subsection{Bounding the $L^2 \to H^s$ Norm of $K_j$}
We will first provide an upper bound for the quantity $\|K_j\|_{L^{2}\to H^s}$ for any $s\in \mathbb{Z}$. To do this, we will first generalize the operator $K_j$, which outputs a scalar-valued function when applied to $f$, to the operator $A_j$, which outputs a vector-valued function instead. Then, using \Cref{lem:linear_map_sobolev_norm_boundary}, we will convert the problem of bounding $\|K_j f\|_{H^s}$, for a function $f$, into the problem of bounding $\|A_j f\|_{L^2}$.
\subsubsection{Generalizing the Operator $K_j$}
Consider a vector-valued function
\[
    \bv(x,y_1): \partial\Omega \times B_a \to \R^d,
\]
where $d$ is the number of different multi-indices $\alpha$ such that $|\alpha| \le |s|$.
We define
\[
    A_j: L^2(\mathcal{B}) \to L^2(\partial\Omega\times\partial\Omega)
\]
such that for function $f$ in~\eref{eq:f def} and $j\geq 1$, we have
\begin{eqnarray}\label{eq:A_operator_formula} \fl
    (A_j f)(x,x_1) = k^{2j}\int_{\mathcal{B}_j}\bv(x,y_1)G(y_1,y_2)\cdots G(y_{j-1},y_j)  \,f\, u_i(y_j,x_1)dy_1\cdots dy_j.
\end{eqnarray}
In other words, the operator $A_j$ is a generalization of $K_j$ obtained by replacing the scalar-valued $G(x,y_1)$ with the vector-valued $\bv(x,y_1)$. We first determine an upper bound on the $L^2$ norm of the operator $A_j$ for $j \ge 1$, and here,
\[
    \|A_j f\|_{L^2}^2 = \sum_{r=1}^m \sum_{r_1=1}^m\int_{\Delta_{r_1}}\int_{\Delta_r}|(A_j f)(x,x_1)|^2dx_r'dx_{r_1}',
\]
where for a pair of local charts $(r,r_1)$, we have $x = (x_r',a_r(x_r'))$ and $x_1 = (x_{r_1}',a_{r_1}(x_{r_1}'))$. We first prove the following lemma which gives an upper bound on $\|A_j\|_{L^2}$. We define
\[
C_a = |B_a|^{\frac{1}{2}}\sup_{y\in B_a}\|u_i(y,\cdot)\|_{L^2(\partial\Omega)}.
\]
\begin{lem}\label{lem:A_L2_L2_bound}
Let $\mu_{A} := k^2\sup_{y\in B_a}\|G(y,\cdot)\|_{L^2(B_a)}$ and 
$\nu_{A} := k^2C_a \sup\limits_{y\in B_a}\|\bv(\cdot,y)\|_{L^2(\partial\Omega)}$. Then, the operator $A_j$ defined by \eref{eq:A_operator_formula} satisfies
\[
    \|A_j\|_{L^2} \le \nu_{A}\,\mu_{A}^{j-1}.
\]
\end{lem}

\begin{proof}

From the Cauchy--Schwarz inequality,
\begin{eqnarray*} \fl
    |(A_j f)(x,x_1)|^2 \le k^{4j}\|f\|_{L^2}^2\int_{\mathcal{B}_j}|\bv(x,y_1)G(y_1,y_2)\cdots G(y_{j-1},y_j)u_i(y_j,x_1)|^2dy_1\cdots dy_j.
\end{eqnarray*}
Let $x = (x_r',a_r(x_r'))$ and $x_1 = (x_{r_1}',a_{r_1}(x_{r_1}'))$. We first bound $\|A_1\|_{L^2}$:
\begin{eqnarray*}
    \|A_1 \eta\|_{L^2}^2 &\le k^4 \|\eta\|_{L^2}^2
    \sum_{r=1}^m \sum_{r_1=1}^m\int_{\Delta_r}\int_{\Delta_{r_1}}\int_{B_a}|\bv(x,y_1)u_i(y_1,x_1)|^2dy_1dx_{r_1}'dx_r'\\
    &= k^4 \|\eta\|_{L^2}^2
\int_{B_a} \|\bv(\cdot, y_1)\|_{L^2(\partial\Omega)}^2\|u_i(y_1,\cdot)\|_{L^2(\partial\Omega)}^2 dy_1\\
    &\le k^2 \|\eta\|_{L^2} C_a \sup_{y_1\in B_a}\|\bv(\cdot,y_1)\|_{L^2(\partial\Omega)}.
\end{eqnarray*}
Thus,  $\|A_1\|_{L^2} \le k^2 C_a \sup\limits_{y_1\in B_a}\|\bv(\cdot,y_1)\|_{L^2(\partial\Omega)}$.
Now, we estimate $\|A_j\|_{L^2}$ where $j \ge 2$:
\begin{eqnarray*} 
    \|A_j f\|_{L^2}^2 &\le& k^{4j}\|f\|_{L^2}^2  \int_{\mathcal{B}_j}|G(y_1,y_2)\cdots G(y_{j-1},y_j)|^2dy_1\cdots dy_j \\
&&  \sup_{y_1\in B_a, y_j\in B_a}\sum_{r=1}^m \sum_{r_1=1}^m\int_{\Delta_r}\int_{\Delta_{r_1}}|\bv(x,y_1)u_i(y_j,x_1)|^2dx_{r_1}'dx_r'.
\end{eqnarray*}
This gives
\[
    \|A_j\|_{L^2} \le k^{2j}\, J_{j-1}\, \sup_{y\in B_a}\|u_i(y,\cdot)\|_{L^2(\partial\Omega)}\sup_{y\in B_a}\|\bv(\cdot,y)\|_{L^2(\partial\Omega)},
\]
where $J_{j-1} := \left(\int_{\mathcal{B}_j}|G(y_1,y_2)\cdots G(y_{j-1},y_j)|^2dy_1\cdots dy_j\right)^{\frac{1}{2}}$. 
As shown in~\cite[Eqn.\ (30)]{moskow2008convergence}, 
\[
    J_{j-1} \le |B_a|^{\frac{1}{2}}\left(\sup_{y\in B_a}\|G(y,\cdot)\|_{L^2(B_a)}\right)^{j-1},
\]
which gives
$\|A_j\|_{L^2} \le \nu_A\, \mu_A^{j-1}$ based on their definitions.
\end{proof}
\begin{rem}
If $d=1$ and $\bv(x,y) = G(x,y)$, then the operators $A_j$ and $K_j$ are the same. Thus, the bound proved in \Cref{lem:A_L2_L2_bound} is a more general version of the bound proved in~\cite{moskow2008convergence}. Our bound is different in that $\nu_A$ contains a factor of $\sup_{y\in B_a}\|u_i(y,\cdot)\|_{L^2(\partial\Omega)}$, which is a minor correction to~\cite[Eqn.\ (25)]{moskow2008convergence}.
\end{rem}

\subsubsection{Bounding $\|K_j\|_{L^2 \to H^s}$} In this section, we will consider the operator 
\[
    K_j: L^2(\mathcal{B}_j) \to H^s(\partial\Omega)\times L^2(\partial\Omega)
\]
defined in~\eref{eq:K_operator_formula} for $s\in\mathbb{Z}$.
Using Lemmas~\ref{lem:linear_map_sobolev_norm_boundary} and~\ref{lem:A_L2_L2_bound}, we will get an upper bound on $\|K_j\|_{L^2 \to H^s}$ for any $s\in\mathbb{Z}$. We define 
\[
    \|K_j f\|_{H^s\times L^2}^2 := \sum_{r_1=1}^m\int_{\Delta_{r_1}}\|(K_j f)(\cdot,x_1)\|_{H^s(\partial\Omega)}^2 dx_{r_1}',
\]
where $x_1 = (x_{r_1}',a_{r_1}(x_{r_1}'))$ and ${\|(K_j f)(\cdot,x_1)\|_{H^s(\partial\Omega)}}$ is defined through \Cref{defn:sobolev_w_space_boundary}.
\begin{lem}\label{lem:L2_W^s2_bound}
Let $\mu_s := k^2\sup_{y\in B_a}\|G(y,\cdot)\|_{L^2(B_a)}$
and $\nu_s := k^2C_a \sup\limits_{y\in B_a}\|G(y,\cdot)\|_{H^s(\partial\Omega)}$ for $s\in \mathbb{Z}$.
Then, the operator $K_j$ in~\eref{eq:K_operator_formula} satisfies $
\|K_j\|_{L^2 \to H^s} \le \nu_s\,\mu_s^{j-1}$.
\end{lem}
\begin{proof}
From \Cref{lem:linear_map_sobolev_norm_boundary}, there exists a linear operator $\mathbf{T}$ acting only on $x$ such that
\[
     \|K_j f\|_{H^s\times L^2} =  \|\mathbf{T}K_j f\|_{L^2(\partial\Omega\times \partial\Omega)}
\]
where $\mathbf{T}K_jf(x,x_1) \in \R^d$ for some positive integer $d$. Then, for fixed $x_1$,
\begin{eqnarray*} 
    \fl
    \mathbf{T}K_jf(\cdot,x_1) &= \mathbf{T}\Big(k^{2j}\int_{\mathcal{B}_j}G(\cdot,y_1)G(y_1,y_2)\cdots G(y_{j-1},y_j)f(y_1,\ldots,y_j)u_i(y_j,x_1)dy_1\cdots dy_j\Big)\\
    \fl
    &= k^{2j}\int_{\mathcal{B}_j}\mathbf{T}G(\cdot,y_1)G(y_1,y_2)\cdots G(y_{j-1},y_j)f(y_1,\ldots,y_j)u_i(y_j,x_1)dy_1\cdots dy_j,
\end{eqnarray*}
where the last equation follows since $\mathbf{T}$ is linear. Letting $\bv(\cdot,y_1) = \mathbf{T}G(\cdot,y_1)$, which is a vector-valued function in $x$, we get
\begin{eqnarray*} \fl
    (A_j f)(x,x_1) = k^{2j}\int_{\mathcal{B}_j}\bv(x,y_1)G(y_1,y_2)\cdots G(y_{j-1},y_j)
    f(y_1,\ldots,y_j)u_i(y_j,x_1)dy_1\cdots dy_j.
\end{eqnarray*}
From \Cref{lem:A_L2_L2_bound}, we have
\[
    \|K_j f\|_{H^s\times L^2} = \|A_j f\|_{L^2} \le \nu_A\,\mu_A^{j-1}\,\|f\|_{L^2},
\]
where $\mu_{A} = \mu_s$, and
\[
    \nu_{A} = k^2C_a \sup_{y\in B_a}\|\mathbf{T}G(\cdot,y)\|_{L^2(\partial\Omega)} \nonumber= k^2C_a \sup_{y\in B_a}\|G(\cdot,y)\|_{H^s(\partial\Omega)} = \nu_s,
\]
following from \Cref{lem:linear_map_sobolev_norm_boundary} and the fact that $G$ is symmetric.
\end{proof}

\subsection{Bounding the $\ha \to \hb$ Norm of $K_j$}
We can now derive a bound on the norm of the operator
\[
    K_j: \ha(\mathcal{B}_j) \to \hb(\partial\Omega)\times L^2(\partial\Omega),
\]
where $\fa$ and $\fb$ are integers and $\fa \ge 0$. 
Consider the function $f$ defined in~\eref{eq:f def} based on $\eta \in \ha_0(B_a)$. Then, based on $\mathcal{B}_j$ defined in~\eref{eq:B def}, we have
\begin{equation*}
    \|f\|_{L^2}^2 = \int_{\mathcal{B}_j} |\eta(y_1)|^2|\eta(y_2)|^2\cdots|\eta(y_j)|^2 dy_1 dy_2 \cdots dy_j = \|\eta\|_{L^2}^{2j}.
\end{equation*}
Because $f$ is the function obtained by copying the same $\eta(y)$ for a total of $j$ times, we can simply define
\begin{equation}\label{def:H_a_norm}
    \|f\|_{\ha} := \|\eta\|_{\ha}^j,
\end{equation}
as done in~\cite{ryan2002introduction}. We now prove the following lemma.
\begin{lem}\label{lem:s_1_to_s_2_operator_bound}
Let $\fa,\fb\in\mathbb{Z}$, and $\fa \ge 0$. Define
\begin{eqnarray} \label{eq:mu_s_1_s_2_def}
    \mu_{\fab} &:= k^2 P(\fa,a,n)\sup_{y\in B_a}\|G(y,\cdot)\|_{L^2(B_a)},\\
    \nu_{\fab} &:= k^2 P(\fa,a,n)C_a  \sup_{y\in B_a}\|G(y,\cdot)\|_{\hb(\partial\Omega)}.\label{eq:nu_s_1_s_2_def} 
\end{eqnarray}
Then, the operator $K_j$ defined by \eref{eq:K_operator_formula} satisfies $\|K_j\|_{\ha\to \hb} \le \nu_{\fab}\,\mu_{\fab}^{j-1}$.
\end{lem}
\begin{proof}
First, note that the following equations hold:
\begin{equation}\label{eq:mu_b nu_b}
    \nu_{\fab} = \nu_\fb P(\fa,a,n), \qquad \mu_{\fab} = \mu_\fb P(\fa,a,n).
\end{equation}
We have $\|f\|_{L^2} = \|\eta\|_{L^2}^j$, and from \eref{def:H_a_norm}, we get $\|f\|_{\ha} = \|\eta\|_{\ha}^j$.
Since $\eta \in \ha_0(B_a)$, we can apply \Cref{lem:applied_poincare_inequality} and obtain $\|\eta\|_{L^2} \le \|\eta\|_{\ha}P(\fa,a,n)$. From \Cref{lem:L2_W^s2_bound}, 
\begin{equation*}
 \|K_j f\|_{\hb\times L^2} \le \|f\|_{L^2}\nu_\fb \mu_\fb^{j-1}\le \|f\|_{\ha}\nu_\fb \mu_\fb^{j-1}P(\fa,a,n)^j= \|f\|_{\ha}\nu_{\fab}\mu_{\fab}^{j-1},
\end{equation*}
and the proof is completed.
\end{proof}
\begin{rem}
Through a method similar to~\cite[Prop.2.1]{moskow2008convergence}, one can obtain an upper bound $1/\mu_{\fab}$ on the radius of convergence of the series in \eref{eq:forward_series} under the $\ha\to \hb$ norm.
\end{rem}

\begin{rem}
From \eref{eq:mu_s_1_s_2_def}-\eref{eq:nu_s_1_s_2_def}, we see that both $\mu_{\fab}$ and $\nu_{\fab}$ decrease as $\fa$ increases while  $\fb$ is fixed. As $\fb$ increases while holding $\fa$ constant, $\mu_{\fab}$ is constant, while $\nu_{\fab}$ increases.
\end{rem}

\section{Inverse Scattering Series}\label{section:inverse_scattering_series}
In this section, we study the properties of the inverse scattering series. In~\sref{sec:classic}, we prove the convergence, stability, and approximation error of the inverse scattering series under the $H^\fb
\rightarrow H^\fa$ operator norm in~\Cref{lem:new_inverse_operator_bound}, \Cref{thm:inverse_series_convergence}, \Cref{thm:inverse_stability}, and \Cref{thm:inverse_approx_error}, using the classic approach in~\cite{moskow2008convergence}. In~\sref{sec:geometric}, we present similar results using the geometric approach in~\cite{hoskins2022analysis}, where the assumptions are less strict. Later in~\sref{sec:discussions}, we will compare analogous theorems under the classic and geometric approaches in detail.

\subsection{Classic Approach}\label{sec:classic}
First, we present theorems and proofs  similar to analogous results in~\cite{moskow2008convergence} under the $L^p\mapsto L^p$ operator norm where $p\geq 2$, which are also derived in generality for Banach spaces in~\cite{bardsley2014restarted}. Moreover, these results have been improved in~\cite{moskow2019inverse}, and their analysis applies to our setting as well. Thus, we omit the proofs of \Cref{thm:inverse_series_convergence} and \Cref{thm:inverse_stability}. However, we provide sharper bounds on the approximation error in \Cref{thm:inverse_approx_error}, so we include its proof and proofs of necessary lemmas in~\ref{sec:appendix} for completeness. This will be important as we analyze how the constants from these theorems change with respect to different choices of $\fa$ and $\fb$ in~\sref{sec:discussions}.

As in~\cite{moskow2008convergence}, we assume that $\eta$ can be expressed as a series based on $\phi = u - u_i$:
\begin{equation}\label{eq:inverse_series_def}
    \eta = \mathcal{K}_1\phi + \mathcal{K}_2\phi\otimes\phi + \mathcal{K}_3\phi\otimes\phi\otimes\phi + \cdots,
\end{equation}
and substituting \eref{eq:original_forward_series} into \eref{eq:inverse_series_def} gives the solution for $\mathcal{K}_j$, as shown in~\cite{moskow2008convergence}:
\begin{eqnarray}
    \mathcal{K}_1 &= K_1^+,  &  \nonumber\\
    \mathcal{K}_j &= -\left(\sum_{m=1}^{j-1}\mathcal{K}_m\sum_{i_1 + \cdots + i_m = j}K_{i_1}\otimes\cdots\otimes K_{i_m}\right)\mathcal{K}_1\otimes\cdots\otimes\mathcal{K}_1, \label{eq:curly_K_j_definition}&
\end{eqnarray}
where $\mathcal{K}_1 = K_1^+$ is the regularized pseudoinverse of the operator $K_1$. We will first analyze the convergence of the inverse scattering series~\eref{eq:inverse_series_def} under the $\hb \to \ha$ norm, where $\fa \ge 0$ and $\fa,\fb$ are integers. Here, $\eta$, which is the perturbation from the homogeneous medium, belongs to $\ha_0(B_a)$ and the scattering data $\phi(x,x_1) = u(x,x_1) - u_i(x,x_1)$ belongs to $\hb(\partial\Omega)\times L^2(\partial\Omega)$ where $x_1$ denotes the source location and $x$ represents the receiver location. We will denote
\[
    \|\mathcal{K}_j\|_{\fba} =  \|\mathcal{K}_j\|_{\hb \to \ha}
\]
for convenience. Next, we will first find an upper bound on $\|\mathcal{K}_j\|_{\fba}$ in~\Cref{lem:new_inverse_operator_bound}.

\begin{lem}\label{lem:new_inverse_operator_bound}
Suppose $j \ge 2$ and $(\mu_{\fab} + \nu_{\fab})\|\mathcal{K}_1\|_{\fba} < 1$. Let
\begin{eqnarray}\label{eq:C_original_definition}
    C := \|\mathcal{K}_1\|_{\fba}\exp\left(\frac{1}{1 - (\mu_{\fab} + \nu_{\fab})\|\mathcal{K}_1\|_{\fba}}\right).
\end{eqnarray}
Then, the operator $\mathcal{K}_j: \left(\hb(\partial\Omega)\times L^2(\partial\Omega)\right)^j \to \ha(B_a)$
defined by \eref{eq:curly_K_j_definition} satisfies
\[
    \|\mathcal{K}_j\|_{\fba} \le C(\mu_{\fab} + \nu_{\fab})^j\|\mathcal{K}_1\|_{\fba}^j. 
\]
Moreover, for all $\phi \in \hb(\partial\Omega)\times L^2(\partial\Omega)$, we have
\[
    \|\mathcal{K}_j \phi\otimes\cdots\otimes\phi\|_{\ha} \le C(\mu_{\fab} + \nu_{\fab})^j\|\mathcal{K}_1 \phi\|_{\ha}^j. 
\]
\end{lem}
We then present \Cref{thm:inverse_series_convergence}, which gives an upper bound on the radius of convergence of the inverse scattering series and estimates the series limit $\widetilde{\eta}$.
\begin{thm}\label{thm:inverse_series_convergence}
Suppose $\|\mathcal{K}_1\|_{\fba} < 1/(\mu_{\fab} + \nu_{\fab})$ and $\|\mathcal{K}_1 \phi\|_{\ha} < 1/(\mu_{\fab} + \nu_{\fab})$.
Then, the inverse scattering series converges under the $\hb \to \ha$ norm. Furthermore, if $N \ge 1$ and $\widetilde{\eta}$ is the limit of the inverse scattering series, the following bound holds:
\begin{eqnarray}\label{eq:inverse_convergence_estimate}
    \left|\left|\widetilde{\eta} - \sum_{j=1}^N \mathcal{K}_j\phi\otimes\cdots\otimes\phi\right|\right|_{\ha} \le C\frac{\left[(\mu_{\fab} + \nu_{\fab})\|\mathcal{K}_1 \phi\|_{\ha}\right]^{N+1}}{1 - (\mu_{\fab} + \nu_{\fab})\|\mathcal{K}_1 \phi\|_{\ha}}\, ,
\end{eqnarray}
where $C$ is defined in \eref{eq:C_original_definition}.
\end{thm}
Next, we present a result bounding the perturbations in $\eta$ based on the perturbations in $\phi$ up to multiplication by an explicit stability constant in~\Cref{thm:inverse_stability}.
\begin{thm}\label{thm:inverse_stability}
Consider scattering data $\phi_1$ and $\phi_2$, and let $M := \max\{\|\phi_1\|_{\hb},\|\phi_2\|_{\hb}\}$.
Let $\eta_1$ and $\eta_2$ denote the corresponding limits of the inverse scattering series, and let
\begin{equation}\label{eq:C_star_definition}
    C^* := \max\left\{\frac{1}{\mu_{\fab} + \nu_{\fab}}, C\right\},
\end{equation}
where $C$ is defined in \eref{eq:C_original_definition}. Furthermore, let
\begin{equation}\label{eq:C_tilde_definition}
    \widetilde{C} := \frac{C^*(\mu_{\fab} + \nu_{\fab})\|\mathcal{K}_1\|_{\fba}}{(1-(\mu_{\fab} + \nu_{\fab})\|\mathcal{K}_1\|_{\fba}M)^2}.
\end{equation}
If $\|\mathcal{K}_1\|_{\fba} < 1/(\mu_{\fab} + \nu_{\fab})$ and $M\|\mathcal{K}_1\|_{\fba} < 1/(\mu_{\fab} + \nu_{\fab})$, the following bound holds:
\[
    \|\eta_1 - \eta_2\|_{\ha} < \widetilde{C}\|\phi_1 - \phi_2\|_{\hb \times L^2}.
\]
\end{thm}

While the inverse scattering series converges under the conditions of \Cref{thm:inverse_series_convergence}, its limit is not equal to $\eta$ in general. In \Cref{thm:inverse_approx_error}, we give an upper bound on the distance between the series limit $\widetilde{\eta}$ and $\eta$.
\begin{thm}\label{thm:inverse_approx_error}
Let $\mathcal{M} := \max\{\|\eta\|_{\ha}, \|\mathcal{K}_1K_1\eta\|_{\ha}\}$ and let
\begin{equation}\label{eq:C_s1_s2_definition}
    C_{\fab} := \frac{C^*(\mu_{\fab} + \nu_{\fab})}{(1 - (\mu_{\fab} + \nu_{\fab})\mathcal{M})^2}.
\end{equation}
If
$\|\mathcal{K}_1\|_{\fba} < 1/(\mu_{\fab} + \nu_{\fab})$, $\|\mathcal{K}_1 \phi\|_{\ha} < 1/(\mu_{\fab} + \nu_{\fab})$, and $\mathcal{M} < 1/(\mu_{\fab} + \nu_{\fab})$, then the approximation error of the partial sum of \eref{eq:inverse_series_def} satisfies the following bound:
\begin{eqnarray*} \fl
    \bigg\|\eta - \sum_{j=1}^N \mathcal{K}_j\phi\otimes\cdots\otimes\phi\bigg\|_{\ha} \le C_{\fab}
    \|(I - \mathcal{K}_1K_1)\eta\|_{\ha} + C\frac{\left[(\mu_{\fab} + \nu_{\fab})\|\mathcal{K}_1 \phi\|_{\ha}\right]^{N+1}}{1 - (\mu_{\fab} + \nu_{\fab})\|\mathcal{K}_1 \phi\|_{\ha}}
\end{eqnarray*}
where $C$ is defined in \eref{eq:C_original_definition}.
\end{thm}
Note that \Cref{thm:inverse_approx_error} also implies the following by taking $N\rightarrow \infty$.
\begin{cor}
Suppose the hypotheses of \Cref{thm:inverse_approx_error} hold. Then, the approximation error of the inverse scattering series satisfies the following bound:
\begin{eqnarray}\label{eq:inverse_approx_error_formula}
    \|\eta - \widetilde{\eta}\|_{\ha} \le C_{\fab}\| \eta - \mathcal{K}_1K_1\eta\|_{\ha},
\end{eqnarray}
where $\widetilde{\eta}$ is the limit of the inverse scattering series and $C_{\fab}$ is defined in \eref{eq:C_s1_s2_definition}.
\end{cor}

\subsection{Geometric Approach}\label{sec:geometric}
Recent work~\cite{hoskins2022analysis} has improved on the results in~\cite{moskow2008convergence,bardsley2014restarted,moskow2019inverse} through a different approach, obtaining estimates on the radius of convergence and approximation error of the inverse scattering series with more relaxed assumptions. The results in~\cite{moskow2008convergence} can easily be adapted to our setup. We thus include statements and omit the proofs which are similar to those in~\cite{moskow2008convergence}.  We will later investigate how these results change with respect to the choices of $\fa$ and $\fb$ and how they are compared with those from the classic approach. 

First, we present \Cref{thm:inverse_converge_geom}, which provides a radius of convergence estimate for the inverse scattering series, with different conditions from \Cref{thm:inverse_series_convergence}. The main difference is that \Cref{thm:inverse_converge_geom} no longer requires an upper bound for $\|\mathcal{K}_1\|_{\fba}$, which is difficult to achieve in realistic settings unless the regularization coefficient is very large~\cite{hoskins2022analysis}.

\begin{thm}\label{thm:inverse_converge_geom}
Let $C_1 = \max\{2,\|\mathcal{K}_1\|_{\fba}\nuab\}$. Then, the inverse scattering series converges if $\|\mathcal{K}_1 \phi\|_{\ha} < r$, where 
\begin{equation}\label{eq:geom_inverse_radius}
    r = \frac{1}{2\muab}\left(\sqrt{16C_1^2 + 1} - 4C_1\right)
\end{equation}
is the radius of convergence.
\end{thm}

Next, we present \Cref{thm:inverse_approx_geom}, which provides an estimate on the approximation error of the inverse scattering series with different conditions than \Cref{thm:inverse_approx_error}.
\begin{thm}\label{thm:inverse_approx_geom}
Suppose that the forward scattering series and the inverse scattering series both converge, so $\|\mathcal{K}_1 \phi\|_{\ha} < r$, where $r$ is defined in~\eref{eq:geom_inverse_radius}. Let $\widetilde{\eta}$ be the sum of the inverse scattering series and $\mathcal{M}_1 := \max\{\|\eta\|_{\ha}, \|\widetilde{\eta}\|_{\ha}\}$. Define the constant $\widetilde{C}_{\fab}$ by 
\begin{equation}\label{eq:C_approx_geom_def}
    \widetilde{C}_{\fab} := \left(1 - \frac{\nuab\|\mathcal{K}_1\|_{\fba}}{(1 + \nuab\|\mathcal{K}_1\|_{\fba} - \muab\mathcal{M}_1)^2} \right)^{-1}.
\end{equation}
If $\mathcal{M}_1 < \frac{1}{\muab}\left(1 - \sqrt{\frac{\nuab\|\mathcal{K}_1\|_{\fba}}{1 + \nuab\|\mathcal{K}_1\|_{\fba}}}\right)$,
then the approximation error can be bounded above as
\begin{eqnarray*} \fl
    \bigg\|\eta - \sum_{j=1}^N \mathcal{K}_j\phi\otimes\cdots\otimes\phi\bigg\|_{\ha} &\le
   \widetilde{C}_{\fab}\|(I - \mathcal{K}_1K_1)\eta\|_{\ha} \\
    &+  \frac{2\muab}{\sqrt{16C_1^2 + 1}} \left(\frac{\|\mathcal{K}_1\phi\|_{\ha}}{r}\right)^{N+1}\frac{1}{1 - \frac{\|\mathcal{K}_1\phi\|_{\ha}}{r}}.
\end{eqnarray*}
\end{thm}

\section{Discussion and Summary of Results}\label{sec:discussions}
In this section, we analyze how the radius of convergence, stability, and approximation error of the inverse scattering series in~\sref{section:inverse_scattering_series} change with respect to $\fa$ and $\fb$. The discussion could shed light on how to choose the metric space for the parameter $\eta$ and the choice of the objective function in the corresponding inverse data-matching problems. In~\sref{sec:diss_classic}, we first interpret the results in~\sref{sec:classic} obtained by the classic approach. In~\sref{sec:compare}, we then compare the results from the classic approach and the geometric approach given in~\sref{sec:geometric}, followed by numerical illustrations on the radius of convergence for the scalar and diffuse waves in~\sref{sec:numerics}. 

\subsection{Discussions on Results in~\sref{sec:classic}}\label{sec:diss_classic}
In this subsection, we analyze and interpret the results presented in~\sref{sec:classic} using the classic approach.
\subsubsection{Radius of Convergence and Assumptions}\label{sec:assumptions}
Based on~\Cref{thm:inverse_series_convergence}, the inverse scattering series will converge if 
\begin{eqnarray*}
    \|\mathcal{K}_1\|_{\fba} < \frac{1}{\mu_{\fab} + \nu_{\fab}}, \qquad \|\mathcal{K}_1\phi\|_{\ha} < \frac{1}{\mu_{\fab} + \nu_{\fab}}.
\end{eqnarray*}
The radius of convergence of the inverse scattering series, $1/ (\mu_{\fab} + \nu_{\fab})$, controls the allowed sizes of $\|\mathcal{K}_1\|_{\fba}$ and $\|\phi\|_{\hb}$ needed to apply our results in~\sref{section:inverse_scattering_series}. An increased radius of convergence is beneficial, as it indicates that the inverse scattering series will converge under a relatively large perturbation in data, $\|\phi\|_{\hb}$. We can thus apply the results from~\sref{section:inverse_scattering_series} to more cases. As $\fa$ increases, both $\mu_{\fab}$ and $\nu_{\fab}$ decrease, causing $1/ (\mu_{\fab} + \nu_{\fab})$ to increase. As $\fb$ decreases, $\nu_{\fab}$ will decrease and $\mu_{\fab}$ remains constant, causing $1/ (\mu_{\fab} + \nu_{\fab})$ to increase. To sum up, the radius of convergence of the \textit{inverse} scattering series increases as $\fa$ increases and $\fb$ decreases. The radius of convergence of the \textit{forward} scattering series increases as $\fa$ increases and does not depend on $\fb$.

To apply the results from \sref{section:inverse_scattering_series}, we will assume that the conditions of \Cref{lem:new_inverse_operator_bound}, \Cref{thm:inverse_series_convergence}, \Cref{thm:inverse_stability}, and \Cref{thm:inverse_approx_error} hold. As the size of $\|\mathcal{K}_1\|_{\fba}$ is determined by regularization~\cite[Remark 3.2]{moskow2008convergence}, we will assume that 
\begin{equation}\label{eq:new_assumption_curly_K_1}
    \|\mathcal{K}_1\|_{\fba} = Q/ (\mu_{\fab} + \nu_{\fab}),
\end{equation}
for a fixed constant $0 < Q < 1$ independent of $\fa$ and $\fb$.

\subsubsection{Analyzing Stability and Approximation Error}\label{sec:results_analysis}
We first analyze how the constant $C$ in \eref{eq:C_original_definition}, and the constant $C^*$ in \eref{eq:C_star_definition}, change with respect to $\fa$, $\fb$ under the above assumptions. Based on the assumption in~\eref{eq:new_assumption_curly_K_1}, we have
\begin{equation}\label{eq:C_mu_and_nu_formula}
    C = \frac{Q \exp \left[(1-Q)^{-1} \right]}{\mu_{\fab} + \nu_{\fab}}, \qquad C^* = \max\left\{\frac{1}{\mu_{\fab} + \nu_{\fab}},\, C \right\},
\end{equation}
implying that both $C$ and $C^*$ are proportional to $1/ (\mu_{\fab} + \nu_{\fab})$, the inverse scattering series radius of convergence. Thus, under the assumption in~\eref{eq:new_assumption_curly_K_1}, we see that $C$ and $C^*$ decrease as $\fa$ decreases and $\fb$ increases.

\paragraph{Stability Constant} 
For scattering data $\phi_1$ and $\phi_2$ and the corresponding limits $\eta_1$ and $\eta_2$ of the inverse scattering series, we show in \Cref{thm:inverse_stability} that if $M = \max\{\|\phi_1\|_{\hb},\|\phi_2\|_{\hb}\}$ and the conditions
\begin{equation}\label{eq:stability_conditions}
    \|\mathcal{K}_1\|_{\fba} < 1/ (\mu_{\fab} + \nu_{\fab}), \qquad M\|\mathcal{K}_1\|_{\fba} <  1/ (\mu_{\fab} + \nu_{\fab})
\end{equation}
hold, then
\[
    \|\eta_1 - \eta_2\|_{\ha} < \widetilde{C}\|\phi_1 - \phi_2\|_{\hb \times L^2},
\]
where $\widetilde{C}$ is defined in \eref{eq:C_tilde_definition}. 
From the assumption in \eref{eq:new_assumption_curly_K_1}, we get
\[
    \widetilde{C} = \frac{C^*Q}{(1-QM)^2}.
\]
Note that $M$ is independent of $\fa$, so $\widetilde{C}$ decreases as $\fa$ decreases. It indicates that we obtain better stability in the inverse scattering problem if we seek the parameter in a weaker (and thus bigger) function space.

However, as $\fb$ increases, $C^*$ decreases while $M$ increases. Furthermore, the rate at which $M$ increases is unknown, as this depends on $\phi_1$ and $\phi_2$, causing the change in $\widetilde{C}$ to be unclear. However, as $\fb$ increases, it is less likely that the second condition in \eref{eq:stability_conditions} holds, in which case the bound on $\widetilde{C}$ may not hold as the inverse scattering series could fail to converge. To this end, it is also beneficial to consider a weaker function space for the data, i.e., to use a weaker norm as the objective function in the corresponding PDE-constrained optimization problem to invert $\eta$ computationally.

\paragraph{Approximation Error} 
Under the assumptions of \Cref{thm:inverse_approx_error}, we provide a bound on the approximation error of the inverse scattering series in \eref{eq:inverse_approx_error_formula}. This bound decreases with the constant $C_{\fab}$, which is defined in \eref{eq:C_s1_s2_definition}.
Based on~\eref{eq:C_mu_and_nu_formula}, we get
\begin{equation}\label{eq:cab}
     C_{\fab} = \frac{\max\left\{1,Q\exp \left[(1-Q)^{-1} \right]\right\}}{(1 - (\mu_{\fab} + \nu_{\fab})\mathcal{M})^2},
\end{equation}
where $\mathcal{M} = \max\{\|\eta\|_{\ha}, \|\mathcal{K}_1K_1\eta\|_{\ha}\}$. Since $Q$ is a constant independent of $\fa$ and $\fb$, $C_{\fab}$ increases if and only if the constant $(\mu_{\fab} + \nu_{\fab})\mathcal{M}$ increases. Next, observe that 
$$
    \|\mathcal{K}_1K_1\eta\|_{\ha} \le \|\mathcal{K}_1\|_{\fba}\cdot \|K_1\|_{\ha\to \hb} \cdot \|\eta\|_{\ha}= \frac{Q\|K_1\|_{\ha\to \hb} \|\eta\|_{\ha}}{\mu_{\fab} + \nu_{\fab}} 
    \le \frac{Q\nu_{\fab}\|\eta\|_{\ha}}{\mu_{\fab} + \nu_{\fab}}, 
$$
in which we applied~\eref{eq:new_assumption_curly_K_1} and \Cref{lem:s_1_to_s_2_operator_bound}. 
Since $Q < 1$, this gives $\|\mathcal{K}_1K_1\eta\|_{\ha} \le \|\eta\|_{\ha}$, implying that $\mathcal{M} = \|\eta\|_{\ha}$. Thus, $(\mu_{\fab} + \nu_{\fab})\mathcal{M}$ increases in $\fb$, since $\|\eta\|_{\ha}$ and $\mu_{\fab}$ are constant in $\fb$ and $\nu_{\fab}$ increases in $\fb$. 

Next, note that $\mu_{\fab} + \nu_{\fab} = (\mu_\fb + \nu_\fb)P(\fa,a,n)$ based on~\eref{eq:mu_b nu_b}. Thus, as $\fa$ changes, $C_{\fab}$ increases if and only if $P(\fa,a,n)\|\eta\|_{\ha}$ increases.
From \Cref{lem:generalized_poincare}, we have 
$$
    \frac{1}{{n-1 \choose n-1}}\|\eta\|_{L^2(B_a)}^2 \le \frac{2a^2}{{n \choose n-1}} \sum_{|\alpha| = 1}\|D^{\alpha}\eta\|_{L^2(B_a)}^2 \le \cdots \le \frac{(2a^2)^{\fa+1}}{{n+\fa \choose n-1}} \sum_{|\alpha| = \fa+1}\|D^{\alpha}\eta\|_{L^2(B_a)}^2.
$$
This implies
$$
    \frac{\sum_{|\alpha| \le \fa} \|D^{\alpha}f\|_{L^2(\Omega)}^2}{\sum_{j = 0}^{\fa} \frac{{n+j-1 \choose n-1}}{(2a^2)^j}} \le \frac{\sum_{|\alpha| \le \fa+1} \|D^{\alpha}f\|_{L^2(\Omega)}^2}{\sum_{j = 0}^{\fa+1} \frac{{n+j-1 \choose n-1}}{(2a^2)^j}}.
$$
Thus, $P(\fa,a,n)\|\eta\|_{\ha}$ is increasing in $\fa$. To sum up, $C_{\fab}$ increases in both $\fa$ and $\fb$.

\subsection{Comparing Classic and Geometric Approaches}\label{sec:compare}
Next, we compare the results from \sref{sec:classic} obtained from a classic approach with similar results from \sref{sec:geometric} obtained using geometric function theory. Certainly, \Cref{thm:inverse_converge_geom} and \Cref{thm:inverse_approx_geom}, obtained through geometric function theory, are improvements in the sense that the assumption on $\|\mathcal{K}_1\|_{\fba}$ is removed. However, for this comparison, we use the same assumptions so that all theorems based on both approaches are valid.

\subsubsection{Radius of Convergence}
We begin by analyzing the radius of convergence of the inverse scattering series. The classic approach, shown in \Cref{thm:inverse_series_convergence}, yields a radius of convergence of $r = 1/(\muab + \nuab)$. On the other hand, the assumption \eref{eq:new_assumption_curly_K_1} implies that $2 > \nuab\|\mathcal{K}_1\|_{\fba}$, which means \Cref{thm:inverse_converge_geom} yields a radius of convergence of $r = (\sqrt{65} - 8)/(2\muab)$. Comparing the two, we see that \Cref{thm:inverse_converge_geom} yields a larger radius of convergence whenever $\nuab > \muab (2\sqrt{65} + 15)$. This occurs when $\fb$ is sufficiently large, based on \eref{eq:mu_s_1_s_2_def} and \eref{eq:nu_s_1_s_2_def}. Based on our numerical setup in~\sref{sec:numerics}, the threshold varies for the particular PDE forward model and occurs when $\fb \approx 1$ for the scalar wave and $\fb \approx 4$ for the diffuse wave, as shown in~\fref{fig:inverse_plots_s2_s1_a}. 
In this figure, we also plot the radius of convergence of the inverse scattering series under our assumptions given by both the classic and the geometric approaches for both diffuse and scalar waves as $\fb$ increases. The radius of convergence given by the geometric approach is constant in $\fb$. Although the $\muab$ values are very different between the diffuse and scalar waves, their geometric radii of convergence are visually close due to the multiplication by the constant $(\sqrt{65} - 8)/2 \approx 0.03$.
\subsubsection{Approximation Error}\label{sec:compare_r}
Next, we compare \Cref{thm:inverse_approx_error} and \Cref{thm:inverse_approx_geom}, which both analyze the approximation error of the inverse scattering series. As they are proven using different approaches, they lead to different bounds and conditions. For this analysis, we assume that the conditions of both theorems hold, and that 
\begin{equation*}
    \mathcal{M}_1 = \mathcal{M} = \|\eta\|_{\ha},\quad {and}\quad
    \mathcal{M} = Q_2/(\muab + \nuab),
\end{equation*}
for a fixed $0 < Q_2 < 1$, similar to the assumption in \eref{eq:new_assumption_curly_K_1}. Under these assumptions, the conditions of \Cref{thm:inverse_approx_error} are satisfied. On the other hand, the condition of \Cref{thm:inverse_approx_geom} can be expressed in terms of $Q$ and $Q_2$ as:
\begin{equation}\label{eq:approx_constraint_geom}\fl
    Q_2 < \frac{\muab + \nuab}{\muab}\left(1 - \sqrt{\frac{\nuab Q}{\muab + \nuab +  \nuab Q}}\right) = \left(1 + R\right)\left(1 - \sqrt{\frac{R Q}{1 + R +  R Q}}\right),
\end{equation}
where $R = \nuab / \muab$, constant in $\fa$ and decreasing as $\fb$ increases. Since we know $0<Q <1$ and $R>0$, the last term in~\eref{eq:approx_constraint_geom} has a global minimum in this region, and the minimum value is around $0.737$.
Furthermore, for sufficiently large $R$, $Q_2$ can be larger than $1$, which is an improvement from~\Cref{thm:inverse_approx_error} which requires $0< Q_2 < 1$. To sum up, \Cref{thm:inverse_approx_error} requires $0 < Q,Q_2 < 1$ while \Cref{thm:inverse_approx_geom} does not have these strict assumptions. On the other hand, \Cref{thm:inverse_approx_geom} poses a condition on the relationship among $R$, $Q$ and $Q_2$; see~\eref{eq:approx_constraint_geom}.

Next, we compare the bounds on the approximation error of the limit of the inverse scattering series. For \Cref{thm:inverse_approx_error}, this bound is controlled by the constant $C_{\fab}$ given in~\eref{eq:cab}, while the bound in \Cref{thm:inverse_approx_geom} is controlled by the constant 
\begin{eqnarray*}\fl
    \widetilde{C}_{\fab} = \left(1 - \frac{\nuab\|\mathcal{K}_1\|_{\fba}}{(1 + \nuab\|\mathcal{K}_1\|_{\fba} - \muab\mathcal{M}_1)^2} \right)^{-1} 
    = \left(1 - \frac{QR(1+R)}{(1 + R + R Q - Q_2)^2} \right)^{-1},
\end{eqnarray*}
under all the assumptions above. As $Q$ approaches $1$, the constant $C_{\fab}$ rapidly blows up. However, there is no such issue with $\widetilde{C}_{\fab}$, provided that the constraint \eref{eq:approx_constraint_geom} is met. Also, under~\eref{eq:approx_constraint_geom}, $\widetilde{C}_{\fab}$ does not blow up even if $Q_2$ approaches $1$. However, $C_{\fab} \rightarrow +\infty$ again in this case. From this perspective, \Cref{thm:inverse_approx_geom} is an improvement over \Cref{thm:inverse_approx_error}. Note that the second error terms of the bound in \Cref{thm:inverse_approx_error} and \Cref{thm:inverse_approx_geom} are controlled by the radius of convergence of the inverse scattering series. The larger the radius, the smaller the error obtained from truncating the series, and as discussed previously, whether an approach leads to an increased radius of convergence depends on $\fb$.

\subsection{Numerical Illustrations}\label{sec:numerics}
In this section, we compare the radius of convergence for both the forward and inverse scattering series in $\R^3$.

\subsubsection{Sobolev Space on a Sphere}
Given $s \in \N$, we define the spherical Sobolev space~\cite{barcelo2020characterization}
\[ \fl
    H^{s}(\mathcal{S}^{n-1}) = \left\{f \in L^2(\mathcal{S}^{n-1}) : \|f\|_{H^{s}(\mathcal{S}^{n-1})}^2 = \|f\|_{L^2(\mathcal{S}^{n-1})}^2 + \|(-\Delta_{\mathcal{S}^{n-1}})^{s/2} f\|_{L^2(\mathcal{S}^{n-1})}^2 < \infty\right\},
\]
where $\Delta_{\mathcal{S}^{n-1}}$ is the Laplace--Beltrami operator on the unit sphere in $\R^n$. We equip the space $H^s$ with the Sobolev norm $\|\cdot\|_{H^s(\mathcal{S}^{n-1})}$. Similar to the definition of the $H^s(\R^n)$ norm through the Fourier transform~\cite{arbogast1999methods}, the spherical Sobolev spaces have an equivalent definition using the spherical harmonic transform, which allows us to extend the definition to the case where $s$ is a real number~\cite{barcelo2021fourier}. More precisely, the space of spherical harmonics of degree $\ell$ on $\mathcal{S}^{n-1}$ has an orthonormal basis $Y_{\ell k}$ for $1 \leq k \leq N(n,\ell)$, where 
\[
N(n,0) = 1, \qquad N(n,\ell) = \frac{(2\ell+n-2)\Gamma(\ell+n-2)}{\Gamma(\ell+1)\Gamma(n-1)}, \quad \ell\geq 1.  
\]
Every $f \in L^2(\mathcal{S}^{n-1})$ can be expressed in a spherical harmonic expansion of the form 
\[
f = \sum_{\ell=0}^\infty \sum_{k=1}^{N(n,\ell)} \hat{f}_{\ell k} Y_{\ell k}, \qquad \hat{f}_{\ell k} = \int_{\mathbb{S}^{n-1}}f \cdot \overline{Y_{\ell k}} \,d\sigma.
\]
The Sobolev space $H^s(\mathcal{S}^{n-1})$ with a real number $s$ is defined by 
\begin{equation}\label{eq:sobolev} 
    H^s(\mathcal{S}^{n-1}) = \bigg\{f\in \mathcal{D}'(\mathcal{S}^{n-1}) : {\|f\|'_{H^s(\mathcal{S}^{n-1})}} = \bigg( \sum_{\ell=0}^\infty \sum_{k=1}^{N(n,\ell)} (1+\ell)^{2s}|\hat{f}_{\ell k}|^2  \bigg)^{\frac{1}{2}}< \infty \bigg\}
\end{equation}
where $\mathcal{D}'(\mathcal{S}^{n-1})$ is the space of distributions on $\mathcal{S}^{n-1}$. We remark that the norms $\|\cdot\|_{H^s(\mathcal{S}^{n-1})}$ and $\|\cdot\|'_{H^s(\mathcal{S}^{n-1})}$ are norm-equivalent but not the same.

\subsubsection{Forward Scattering Series}
\begin{figure}[tpb]
\centering
\subfloat[]
{\includegraphics[width=0.49\textwidth]{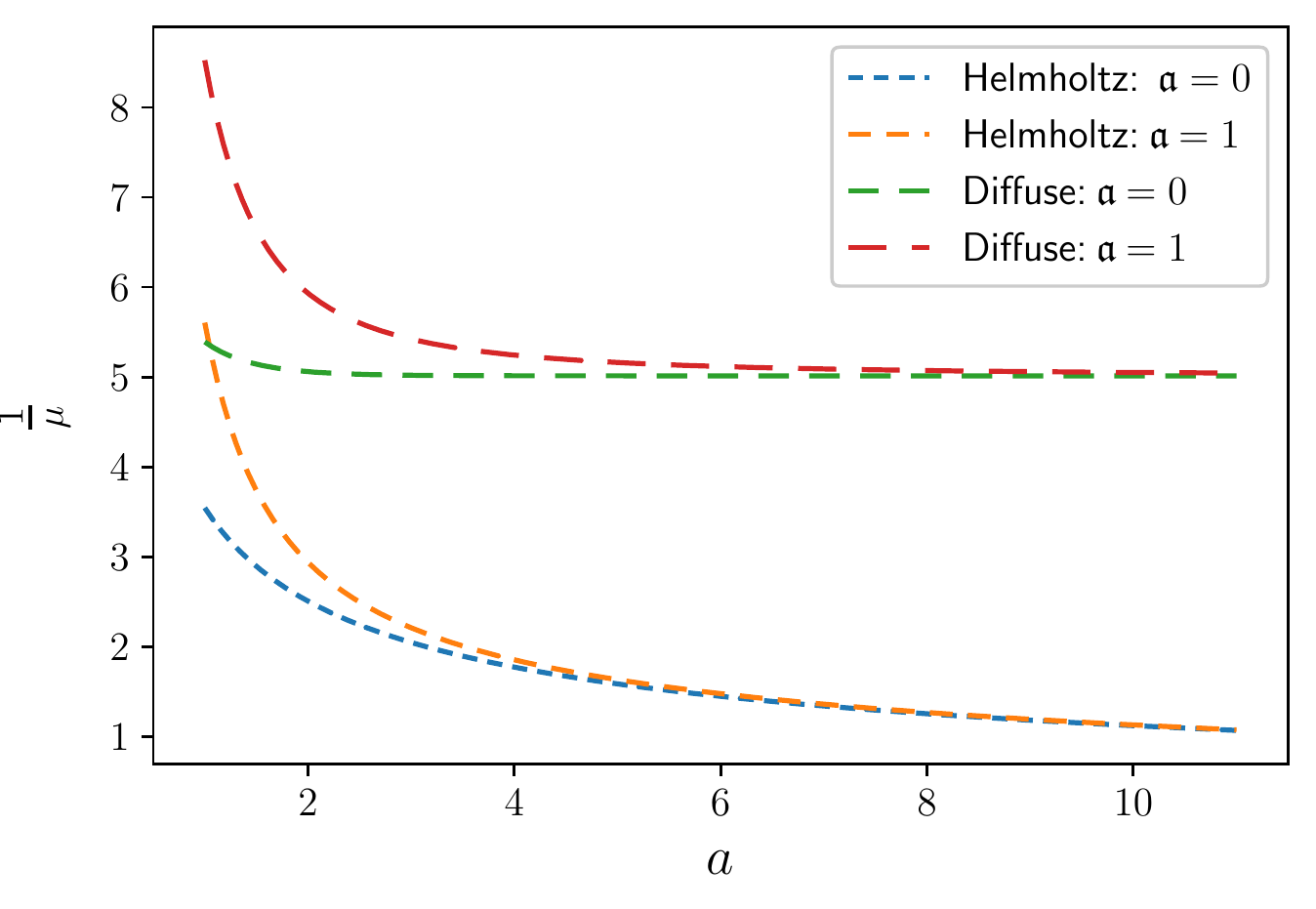}\label{fig:forward_comparison_plots_a}}
\subfloat[]{\includegraphics[width=0.49\textwidth]{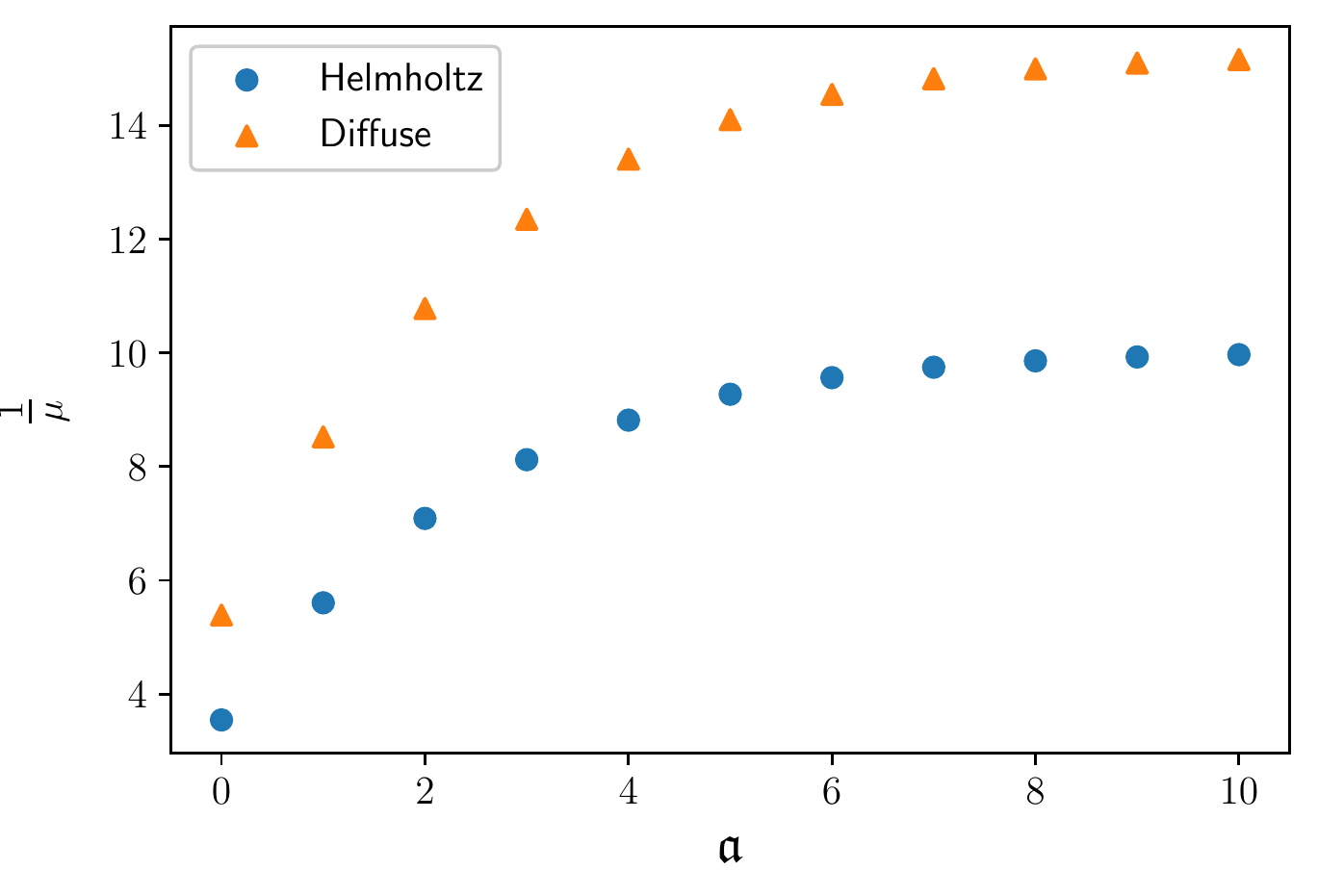}\label{fig:forward_comparison_plots_b} }
\caption{The radius of convergence of the forward scattering series (a)~for Helmholtz and diffuse wave equations under the $\ha \to \hb$ norm as $a$ (the radius of ball $B_a$) ranges from $1$ to $11$, and (b) as $\fa$ ranges from $0$ to $10$.} 
\label{fig:forward_comparison_plots}
\end{figure}
We compute $\mu_{\fab}$ in $\R^3$ for varying $a$ (the radius of $B_a$) and $\fa$ (which controls the parameter space). First, recall that 
\[
    \mu_{\fab} = k^2 P(\fa,a,n)\sup_{y\in B_a}\|G(y,\cdot)\|_{L^2(B_a)}
\]
and note that $\mu_{\fab}$ does not depend on $\fb$. 
In the case of the Helmholtz equation~\eref{eq:helm-2}, we have $G(x,y) = \exp(ik|x-y|)/(4\pi |x-y|)$. Then, we have
\begin{eqnarray*}
    \mu_{\fab}
    &= k^2 P(\fa,a,n)\sup_{y\in B_a}\left(\int_{B_a}\frac{1}{16\pi^2 |x-y|^2} dx\right)^{\frac{1}{2}}&\\
    &= \frac{k^2}{4\pi} P(\fa,a,n)\left(\int_{B_a}\frac{1}{|x|^2} dx\right)^{\frac{1}{2}} 
    = k^2 P(\fa,a,n)\left(\frac{a}{4\pi}\right)^{\frac{1}{2}}.
\end{eqnarray*}
For the diffuse wave~\eref{eq:diffuse eqn}, we have $G(x,y) = \exp(-k|x-y|)/(4\pi |x-y|)$, which gives
\begin{eqnarray*} 
\mu_{\fab}
    &= k^2 P(\fa,a,n)\sup_{y\in B_a}\left(\int_{B_a}\frac{\exp(-2k|x-y|)}{16\pi^2 |x-y|^2} dx\right)^{\frac{1}{2}}\\
    &= \frac{k^2}{4\pi} P(\fa,a,n)\left(\int_{B_a}\frac{\exp(-2k|x|)}{|x|^2} dx\right)^{\frac{1}{2}}\\
    &= k^2 \exp\left(-\frac{ka}{2}\right) P(\fa,a,n) \left(\frac{\sinh(ka)}{4\pi k}\right)^{\frac{1}{2}}.
\end{eqnarray*}
Again, we can see that $\mu_{\fab}$ does not depend on $\hb$ for both equations.

We compare the values of $1/\mu_{\fab}$, the radius of convergence of the forward scattering series, for the Helmholtz equation and the diffuse wave in \fref{fig:forward_comparison_plots}, where we set $k=1$. 
First, we let $a$ increase, with $\fa \in \{0,1\}$ fixed.  As shown in~\fref{fig:forward_comparison_plots_a}, while the radius of convergence for the diffuse wave equation is bounded below as $a$ increases, the radius of convergence for the Helmholtz equation goes to $0$ as $a$ increases. 
Also, note that for both the Helmholtz and diffuse wave equations, $1/\mu_{\fab}$ increases as $\fa$ increases, which aligns with our analysis in~\sref{sec:results_analysis}. The amount of increase becomes more significant as $a$ becomes smaller (i.e., the support of the parameter perturbation $\eta$ becomes smaller), as a result that $P(\fa,a,n)$ defined in~\eref{eq:psan} decreases in $a$. 

Next, in~\fref{fig:forward_comparison_plots_b}, we plot the radius of convergence for $\fa \in \{0,1,2,\ldots,10\}$ for both the Helmholtz and diffuse wave equations.  We set $a = 1$ and $k=1$. Since $P(\fa,a,n)$ increases in $\fa$, we see that $1/\mu_{\fab}$ increases in $\fa$, which is also shown in~\fref{fig:forward_comparison_plots_b}. 
The value of $1/\mu_{\fab}$ appears to converge for both the Helmholtz and diffuse wave equations as $\fa$ increases when $a = 1$. For $a \le 1/\sqrt{2}$, the series in~\eref{eq:psan} does not converge, so for such $a$, the value of $1/\mu_{\fab}$ would grow arbitrarily large as $\fa$ increases. The specific constant $1/\sqrt{2}$ is due to the coefficient of $2$ in \Cref{thm:poincare_inequality}, which is not sharp and can be improved with more precise estimates of the Poincar\'e constant.

\subsubsection{Inverse Scattering Series}

\begin{figure}[tpb]
\centering
\subfloat[]{\includegraphics[width=0.32\textwidth]{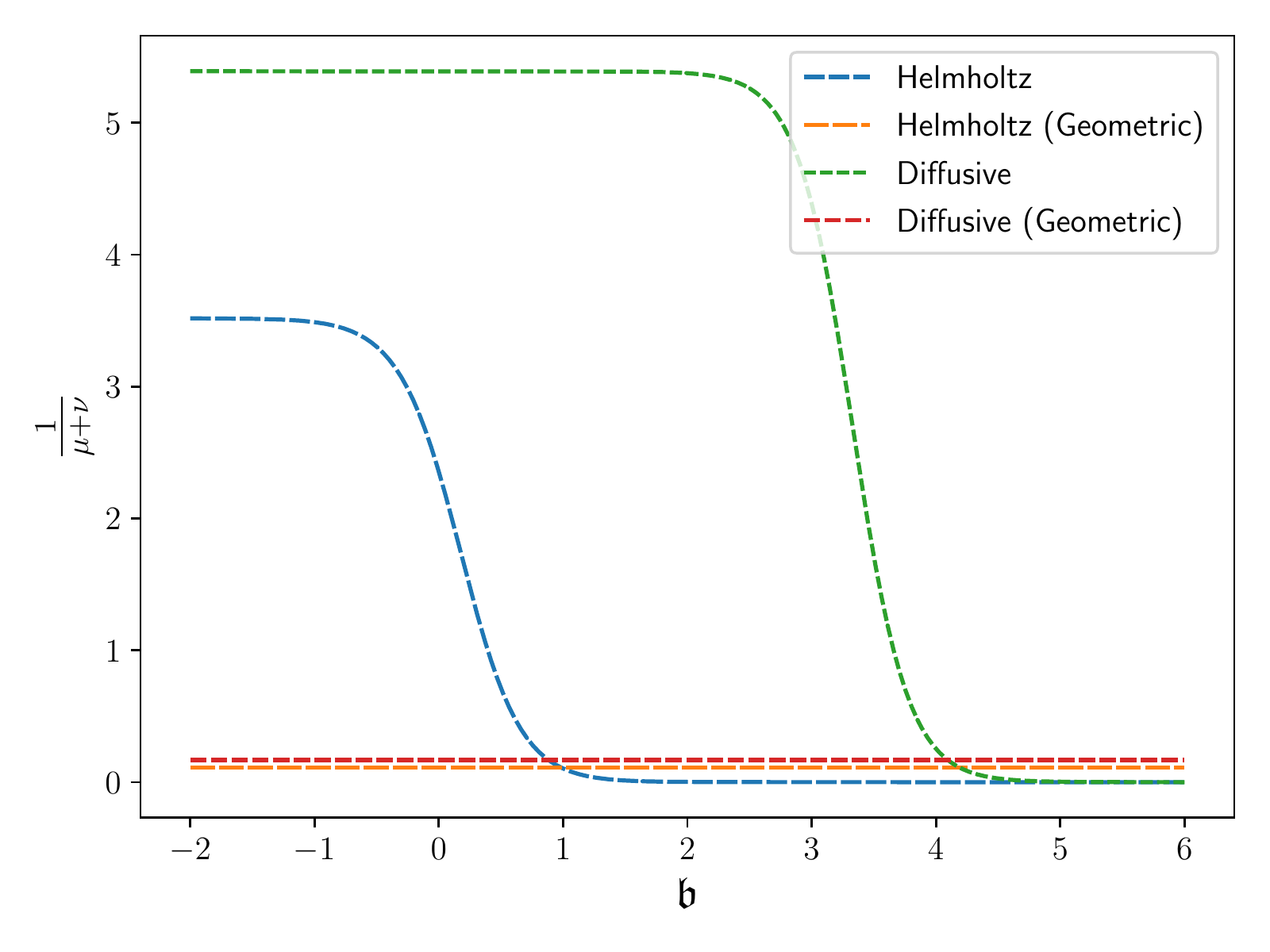}\label{fig:inverse_plots_s2_s1_a}}
\subfloat[]{\includegraphics[width=0.32\textwidth]{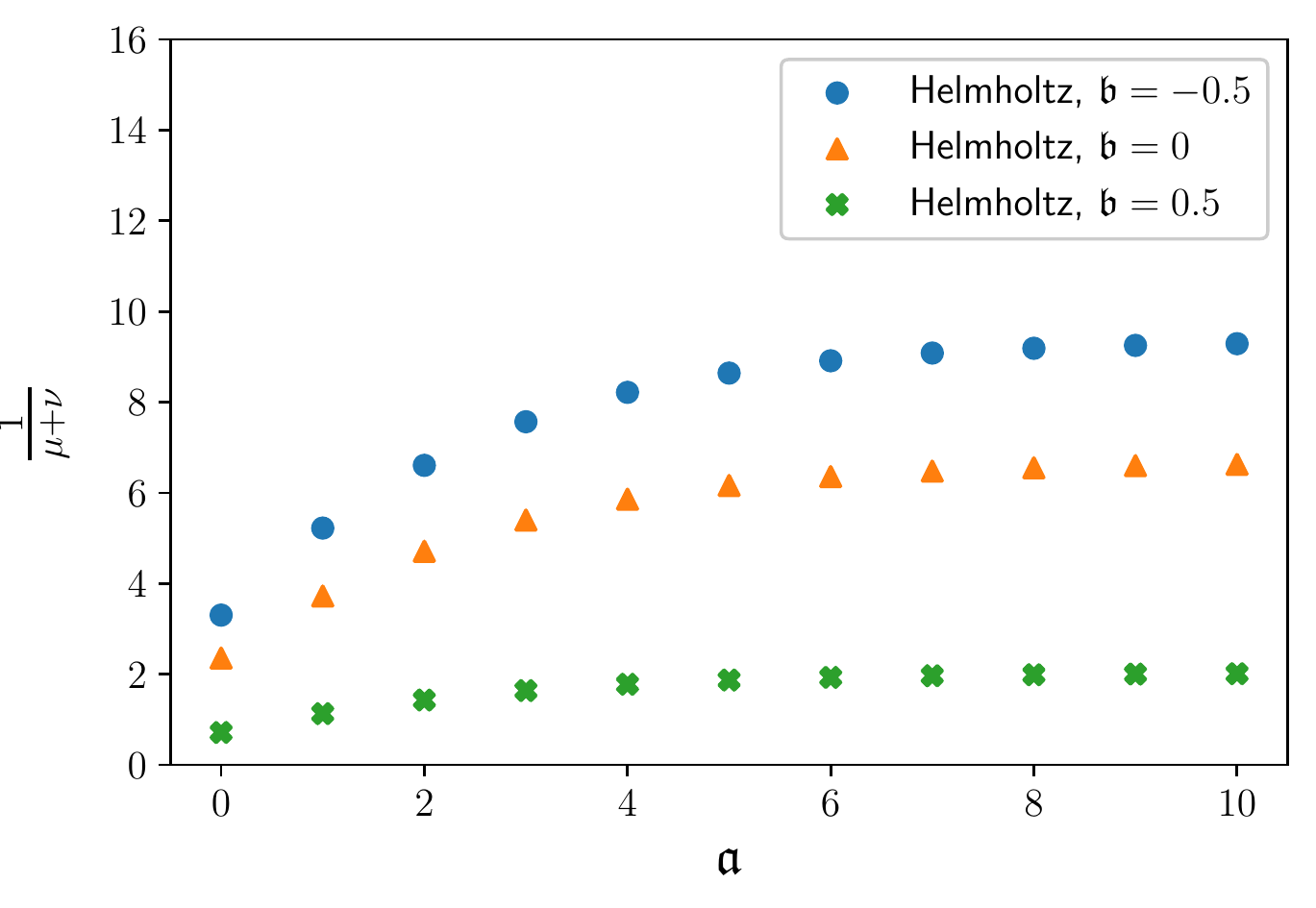}\label{fig:inverse_plots_s2_s1_b}}
\subfloat[]{\includegraphics[width=0.32\textwidth]{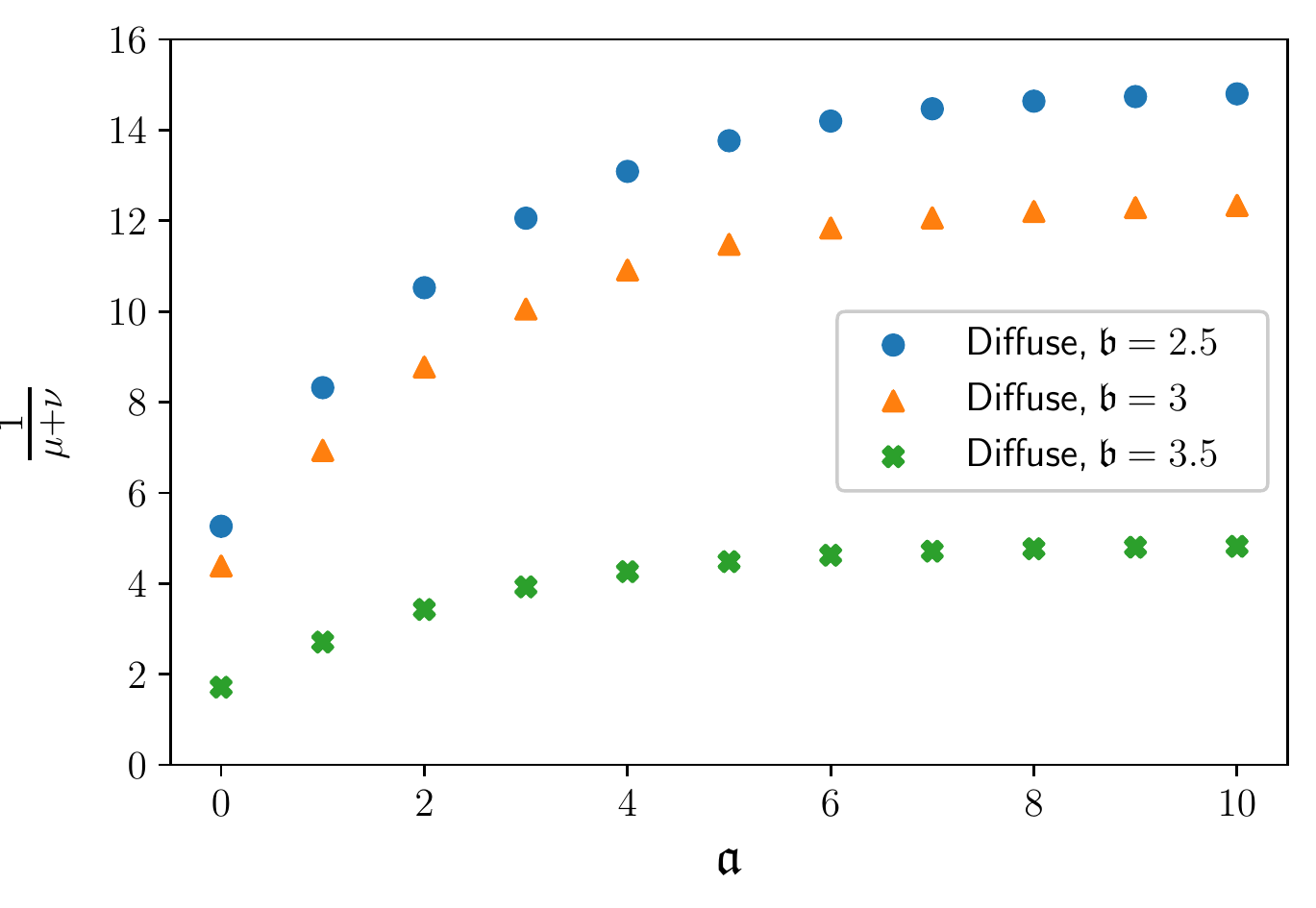}\label{fig:inverse_plots_s2_s1_c}}
\caption{The radius of convergence of the inverse scattering series, (a)~as $\fb$ ranges from $-2$ to $5$, (b)~for the Helmholtz equation as $\fa$ ranges from $0$ to $10$ and (c)~for the diffuse wave equation as $\fa$ ranges from $0$ to $10$. The radius of convergence obtained using the geometric approach is also shown in (a), under the assumption that the conditions of both~\Cref{thm:inverse_series_convergence} and \Cref{thm:inverse_converge_geom} hold. We set $a=1$ and $k=1$ in all three plots.}
\label{fig:inverse_plots_s2_s1}
\end{figure}
We compute $\mu_{\fab}$ and $\nu_{\fab}$ assuming that $\Omega$ is a ball of radius $100$ centered at zero and $B_a$ is centered at $(98,0,0)$ with radius $a=1$ unless otherwise specified. We also assume that $u_i(x,x_1) = G(x,x_1)$ by setting the source to be a delta function. Since $u_i = G$, we get 
\begin{equation*} 
    \nu_{\fab} = k^2 P(\fa,a,n)|B_a|^{\frac{1}{2}}\left(\sup_{y\in B_a}\|G(y,\cdot)\|_{L^2(\partial\Omega)}\right)\left(\sup_{y\in B_a}\|G(y,\cdot)\|_{\hb(\partial\Omega)}\right).
\end{equation*}
Since we are considering the case where $\Omega$ is a ball, we use the definition from \eref{eq:sobolev} to compute $\|G(y,\cdot)\|_{L^2(\partial\Omega)}$ and $\|G(y,\cdot)\|_{\hb(\partial\Omega)}$ for any $y \in B_a$ and $\fb\in\mathbb{R}$. We use the software SHTools~\cite{wieczorek2018shtools} to perform spherical transforms.
The radius of convergence of the inverse scattering series $(\mu_{\fab} + \nu_{\fab})^{-1}$, for both the Helmholtz and diffuse wave equations, are plotted in~\fref{fig:inverse_plots_s2_s1}. We set $a=1$ and $k=1$. In~\fref{fig:inverse_plots_s2_s1_a}, we set $\fa = 0$ and compute $(\mu_{\fab} + \nu_{\fab})^{-1}$ as $\fb$ ranges from $-2$ to $6$. Due to the factor $\|G(y,\cdot)\|_{\hb(\partial\Omega)}$, which grows exponentially based on~\eref{eq:sobolev}, $\nu_{\fab}$ increases exponentially in $\fb$. Thus, the radius of convergence for the inverse scattering series decreases exponentially in $\fb$ for both the Helmholtz and diffuse wave equations. Furthermore, since $(\mu_{\fab} + \nu_{\fab})^{-1}$ is bounded above by $1/\mu_{\fab}$, and the latter is constant in $\fb$, it explains why the decay of  $(\mu_{\fab} + \nu_{\fab})^{-1}$ is similar to that of a sigmoid function. In~\fref{fig:inverse_plots_s2_s1_a}, we also plot the radius of convergence obtained using the geometric approach in~\Cref{thm:inverse_converge_geom}, which is constant in $\fb$ under the same assumptions of \Cref{thm:inverse_series_convergence}.

Also, the radius of convergence for the inverse scattering series is much less sensitive to increases in $\fb$ for the diffuse wave equation than in the case of the Helmholtz equation. This is due to the exponential decay of the Green's function for the diffuse wave equation, as shown in \eref{eq:greens_diffuse}. Since the radius of the domain is relatively large, the terms $\|G(y,\cdot)\|_{L^2(\partial\Omega)}$ and $\|G(y,\cdot)\|_{\hb(\partial\Omega)}$ are small due to this exponential decay. This causes $\nu_{\fab}$ to be far smaller, only having a visible effect on the radius of convergence when $\fb > 2$. We do not see much difference in the radius obtained from the geometric approach due to the multiplication by a small constant; see~\sref{sec:compare_r} for details.

In figures~\ref{fig:inverse_plots_s2_s1_b} and~\ref{fig:inverse_plots_s2_s1_c}, we fix $\fb$ and compute the radius of convergence of the inverse scattering series as $\fa$ increases from $0$ to $10$. 
Note that using different choices of $\fb$ causes the radius of convergence for the Helmholtz equation to change significantly. 
On the other hand, for the diffuse wave equation, $\fb$ can be safely increased to $2.5$ without significantly decreasing the radius of convergence. 
Due to the factor of $P(\fa,a,n)$ in both $\mu_{\fab}$ and $\nu_{\fab}$, the radius of convergence increases as $\fa$ increases, similar to~\fref{fig:forward_comparison_plots_b}. 
Similarly, the behavior of the radius of convergence as $\fa$ grows large depends on $a$. When $a \le 1/\sqrt{2}$, $P(\fa,a,n)$ is unbounded as $\fa$ increases.

\begin{figure}[tpb]
\centering
\subfloat[]{\includegraphics[width=0.33\textwidth]{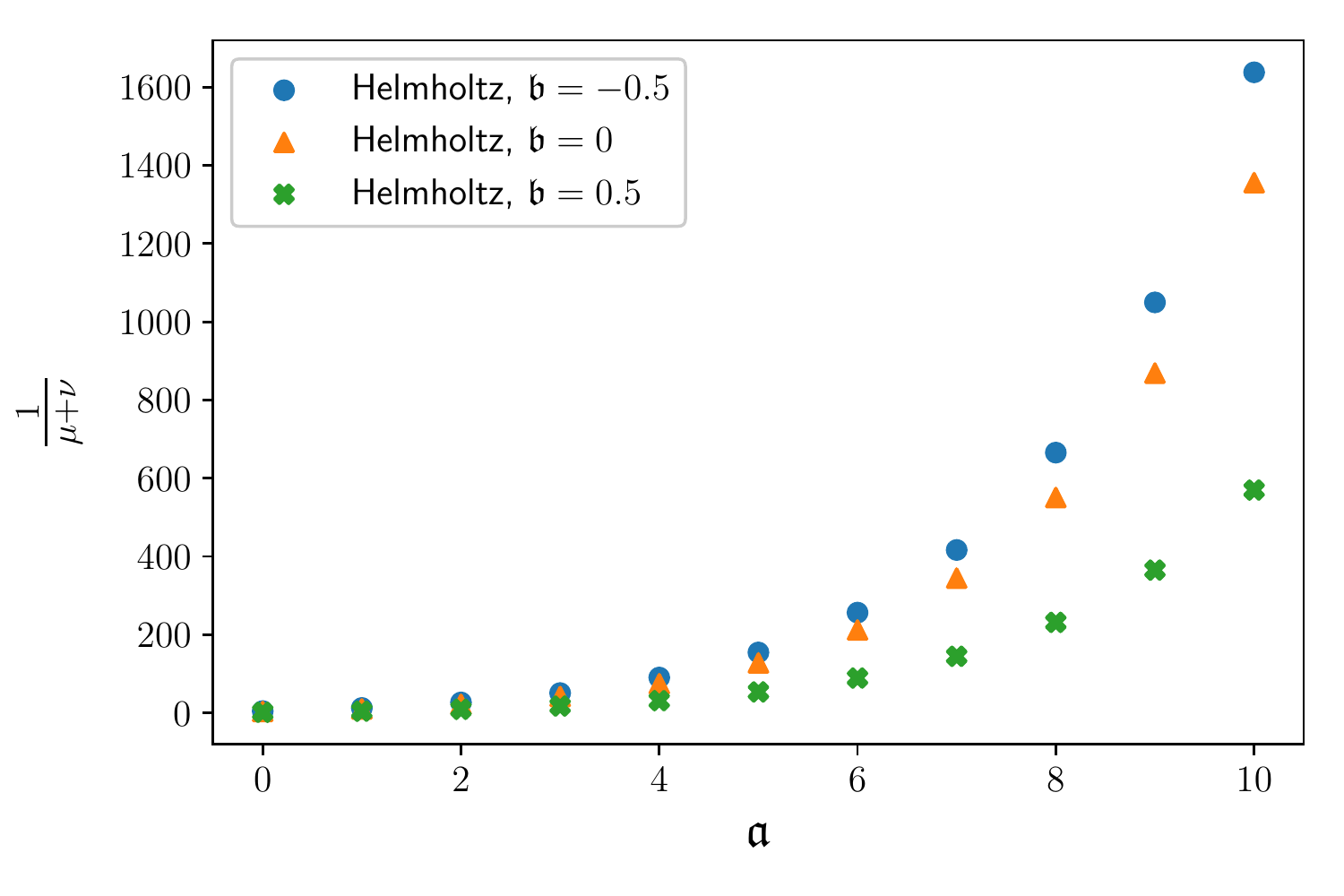}\label{fig:inverse_plots_a1}}
\subfloat[]{\includegraphics[width=0.32\textwidth]{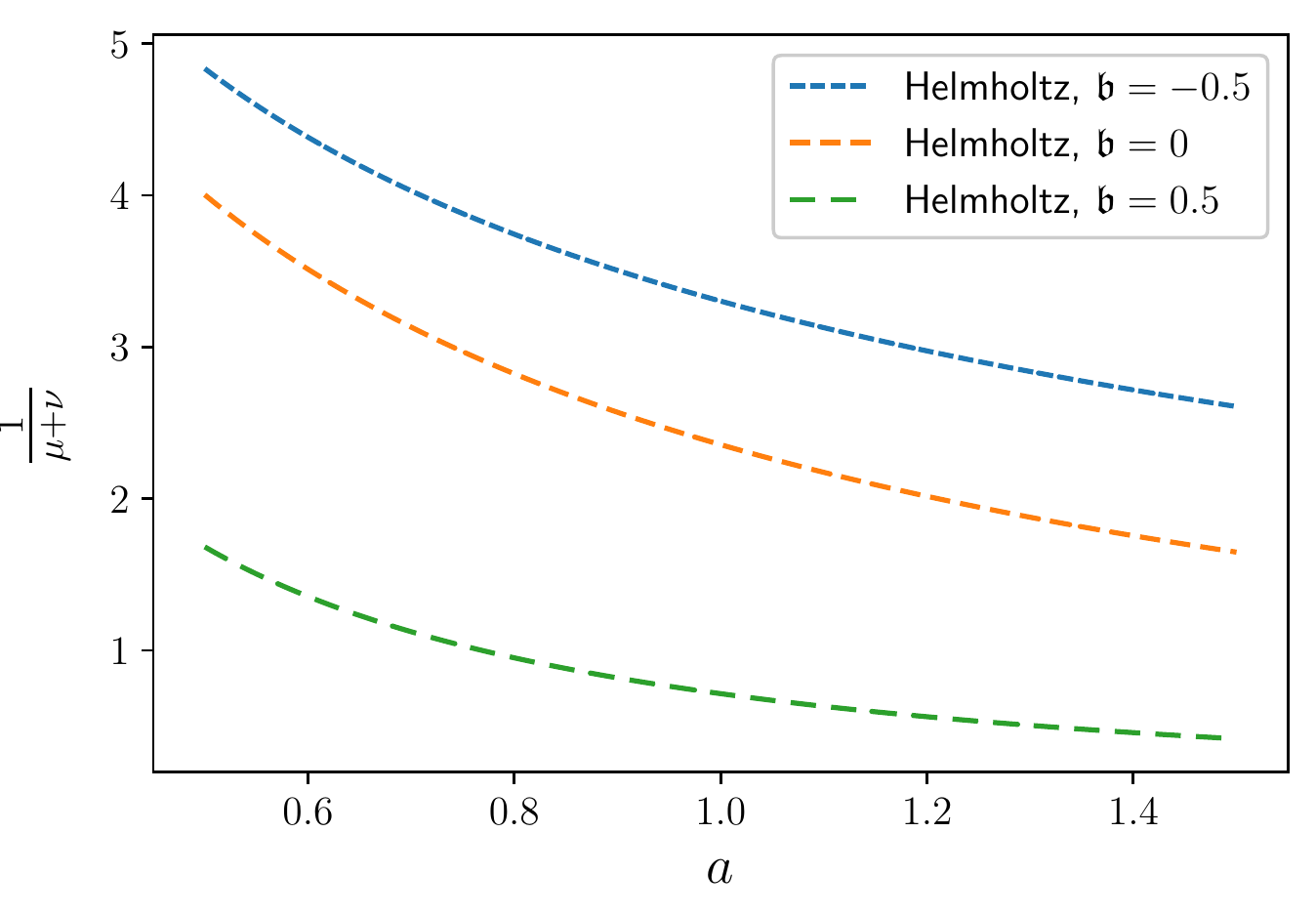}\label{fig:inverse_plots_a2}}
\subfloat[]{\includegraphics[width=0.32\textwidth]{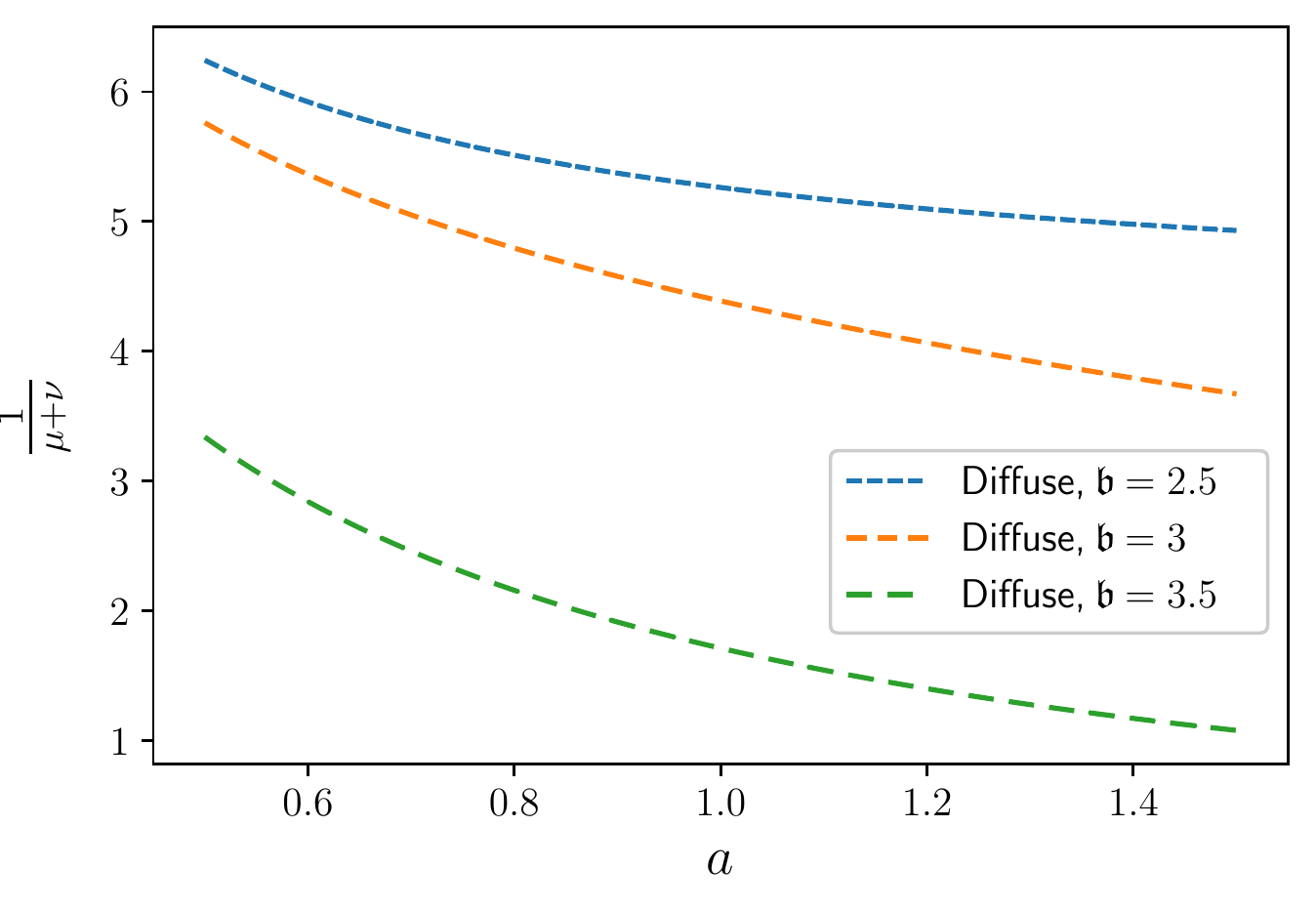}\label{fig:inverse_plots_a3}}
\caption{The radius of convergence of the inverse scattering series (a) as $\fa$ ranges from $0$ to $10$ when $a = 0.5$, (b) for the Helmholtz equation as $a$ ranges from $0.5$ to $1.5$, and (c) for the diffuse wave equation as $a$ ranges from $0.5$ to $1.5$. We set $k=1$ in all plots.}
\end{figure}

In~\fref{fig:inverse_plots_a1}, we demonstrate how the radius of convergence behaves as $\fa$ increases for smaller $a$ (we set $a = 0.5$). Since $P(\fa,a,n)$ is unbounded as $\fa$ increases if $a$ is sufficiently small, the radius of convergence of the inverse scattering series will grow arbitrarily large as $\fa$ increases. In figures~\ref{fig:inverse_plots_a2} and~\ref{fig:inverse_plots_a3}, we plot the radius of convergence for both the Helmholtz and diffuse wave equation as $a$ ranges from $0.5$ to $1.5$. We set $\fa = 0$ since the only effect of $\fa$ on $\mu_{\fab}$ and $\nu_{\fab}$ is an extra multiplicative factor of $P(\fa,a,n)$. The radius of convergence decreases as $a$ increases, and while changing $\fb$ has a significant effect on the Helmholtz equation, $\fb$ needs to be at least $2.5$ to impact the radius of convergence in the case of the diffuse wave equation.

\section{Numerical Inversion Examples}\label{sec:numerical_inversion}
In this section, we present some numerical examples to show the impact of choosing different $H^\fa$ for the parameter space and $H^\fb$ for the data space in an optimization framework to solve the inverse scattering problems.

\subsection{The Gradient Formulation}
To begin with, we explain first how the chosen pair $(\fa,\fb)$ changes the gradient calculation in the PDE-constrained optimization framework. Consider the optimization problem
\begin{equation}\label{eq:obj}
   \min_{\eta \in H^{\fa}} J_\fb(\eta) = \frac{1}{2}\|\mathcal{F}(\eta) - \phi\|_{H^\fb}^2,
\end{equation}
where $\phi$ is the observed data, $\mathcal{F}$ is the forward operator that maps the parameter $\eta$ to the data, which can be represented in the form of~\eref{eq:forward_series} and~\eref{eq:original_forward_series}. We denote by $\mathcal{G}_{\fab}(\eta)$ the gradient of the objective function $J_{\fb}(\eta)$.


By the Riesz representation theorem, $\forall \eta'$, we have
\begin{equation*}\centering
(\mathcal{G}_{\fab}(\eta), \eta' )_{H^{\fa}} = \lim_{\epsilon\rightarrow 0} \frac{J_\fb(\eta + \epsilon \eta') - J_\fb(\eta)}{\epsilon} = \int_\Omega \frac{\delta J_\fb}{\delta \eta} \eta' dx = (\mathcal{G}_{0\fb}(\eta), \eta' )_{L^2},
\end{equation*}
where $ \mathcal{G}_{0\fb}(\eta) = \frac{\delta J_\fb}{\delta \eta}$ denotes the gradient assuming $\eta \in L^2$. We then have the relation
\[
\mathcal{G}_{\fab}(\eta) = (I - \Delta)^{-\fa} \,  \mathcal{G}_{0\fb}(\eta).
\]
As a result, for the same objective function, the gradient with the assumption that $\eta \in H^1$ (i.e., $\eta$ has a higher regularity) is smoother than assuming that $\eta \in L^2$. Similarly, the $L^2$ gradient $\frac{\delta J_\fb}{\delta \eta}$ depends on the choice of $\fb$. Recall that $K_1$ is the linearization of the nonlinear forward operator $\mathcal{F}$. We then have
\[
\frac{\delta J_\fb}{\delta \eta} = K_1^* (I - \Delta)^{\fb} \left(\mathcal{F}(\eta) - \phi \right),
\]
where $K_1^*$ denotes the adjoint operator of $K_1$. If we use $\mathcal{K}_1$ instead of $K_1^*$, it yields the Gauss--Newton method. Combining the above equations, we have 
\begin{equation}\label{eq:Grad_ab}
    \mathcal{G}_{\fab}(\eta) = (I - \Delta)^{-\fa} K_1^* (I - \Delta)^{\fb} \left(\mathcal{F}(\eta) - \phi \right).
\end{equation}
We remark that the two Laplacian operators $\Delta$ in~\eref{eq:Grad_ab} above do not act on the same domain. The first one acts on the model parameter space, while the second one is defined over the data space. Nevertheless, it is evident from~\eref{eq:Grad_ab} that both $\fa$ and $\fb$ directly change the gradient in the PDE-constrained optimization. 

Similarly, to calculate the linear action of $\mathcal{K}_1$, the pseudo-inverse of $K_1$ defined in~\eref{eq:curly_K_j_definition}, it is equivalent to another optimization problem
\begin{equation}\label{eq:obj_linear}
   \min_{\eta \in H^{\fa}}  \frac{1}{2}\|K_1 \eta - \phi\|_{H^\fb}^2.
\end{equation}
Again, the solution changes with respect to the chosen $(\fa,\fb)$, especially when the data is noisy. After discretization, the two norms are reflected in two weight matrices in a least-squares problem. We refer to~\cite{ren2022generalized} where the impact of the weights in the resulting discrete problem is rigorously analyzed under different assumptions on $K_1$.

\begin{figure}
\centering
\subfloat[The true $\eta$]
{\includegraphics[width=0.24\textwidth]{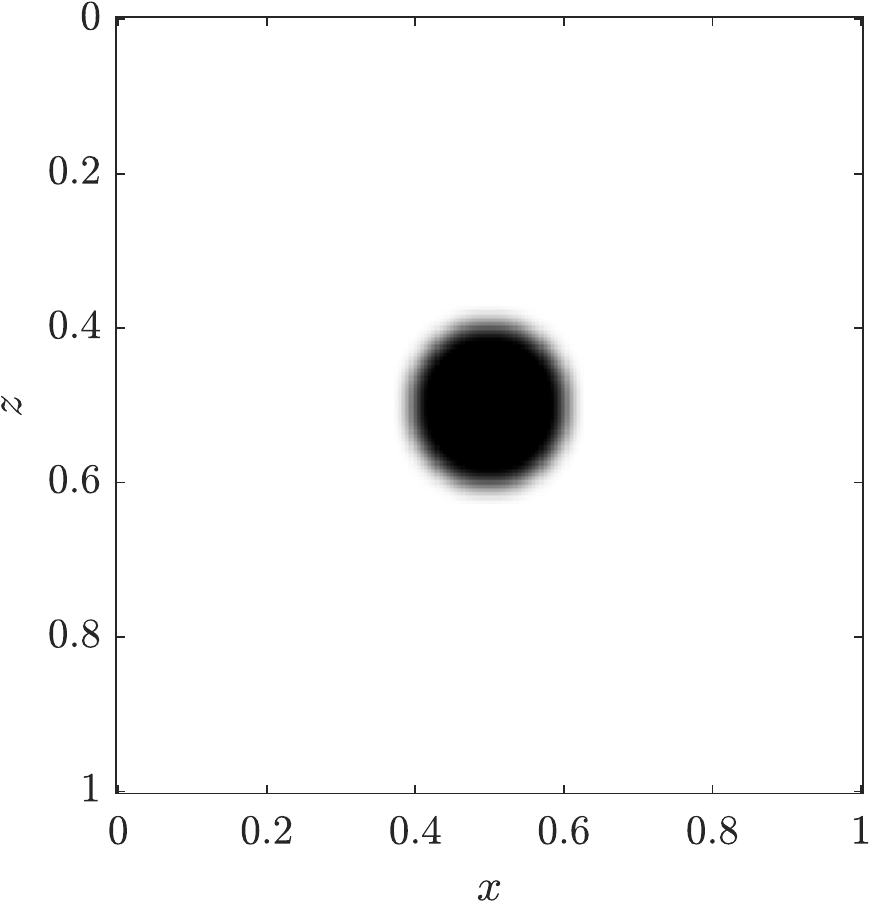}\label{fig:bad true}}
\subfloat[$\fa = \fb = 0$]
{\includegraphics[width=0.24\textwidth]{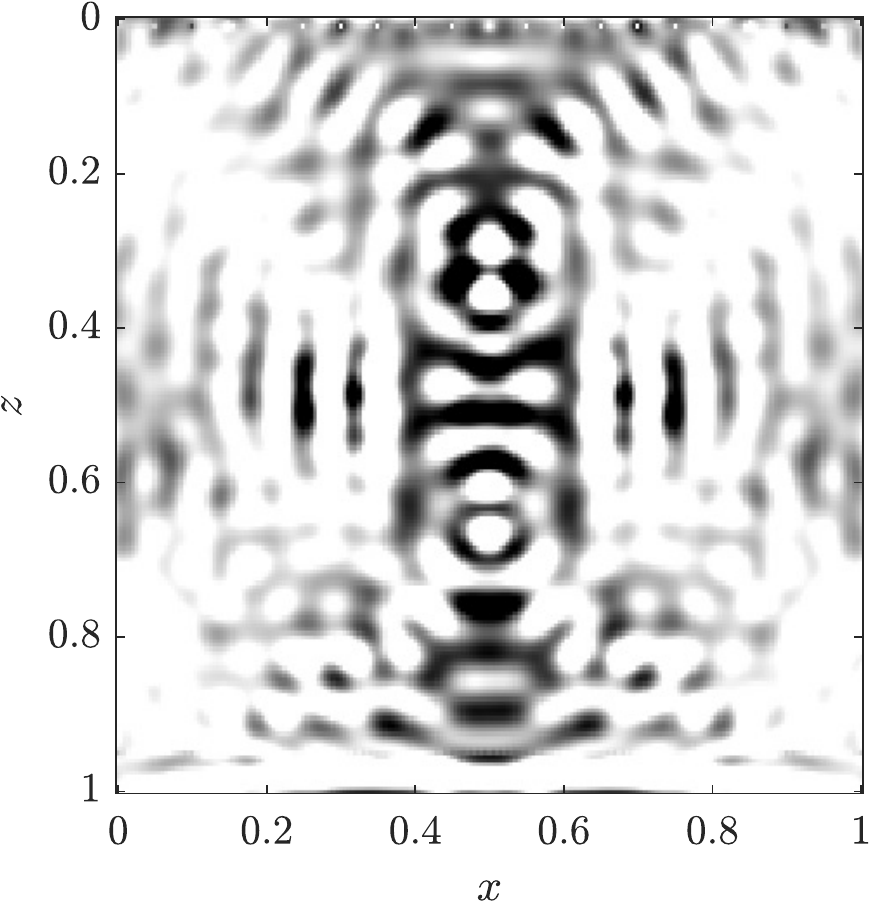}\label{fig:bad a0}}
\subfloat[$\fa = 1$, $\fb = 0$]{\includegraphics[width=0.24\textwidth]{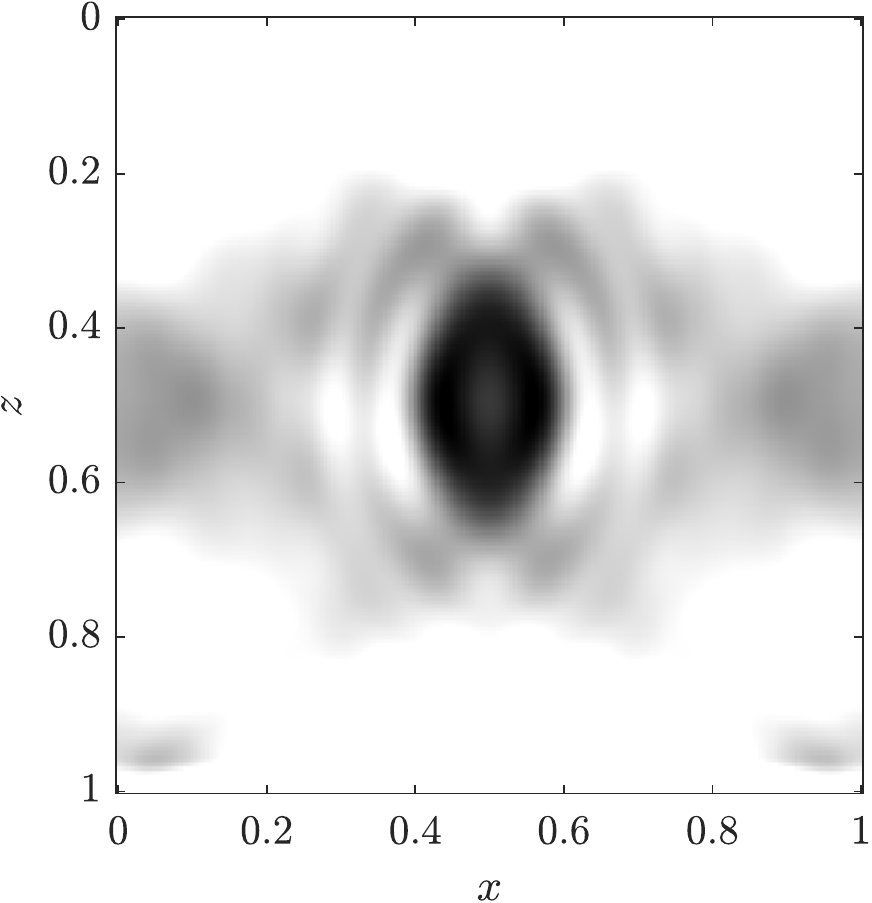}\label{fig:bad a = 1}}
\subfloat[$\fa = 0$, $\fb = -1$]{\includegraphics[width=0.24\textwidth]{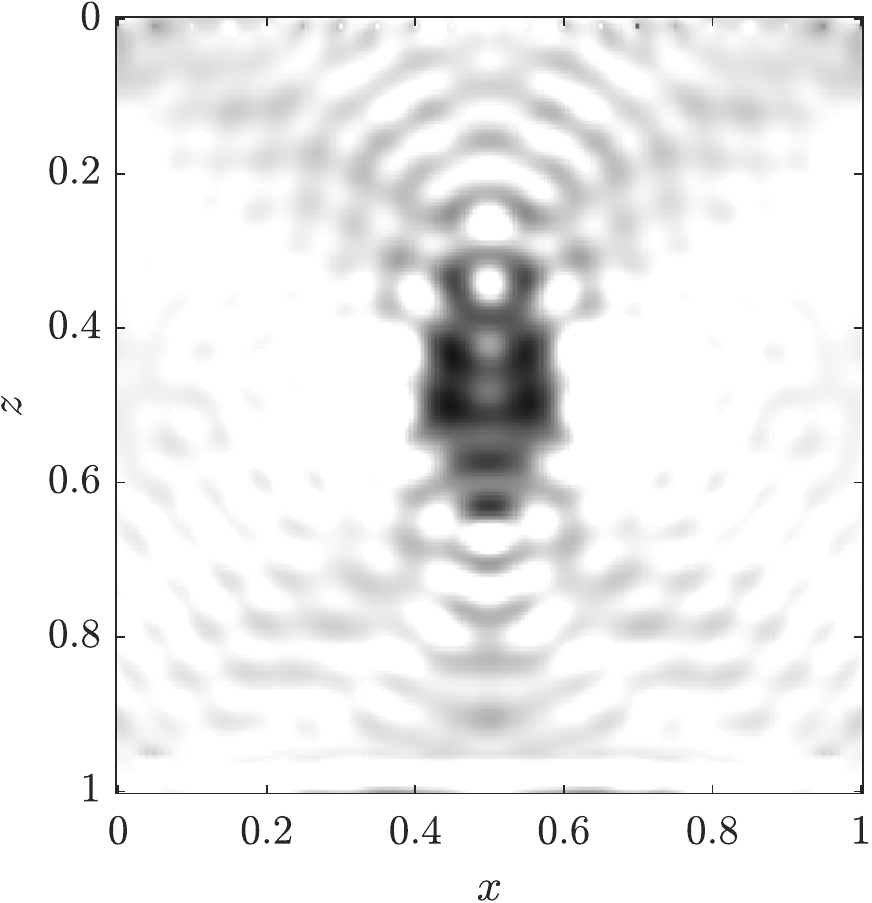}\label{fig:bad b = -1}}
\caption{(Setting One) (a): The true scatterer $\eta$. (b)-(d): The reconstructed results after $100$ iterations of L-BFGS algorithms under different $(\fa,\fb)$ choices.}
\label{fig:bad recovery}
\end{figure}

\begin{figure}
\centering
\subfloat[The true $\eta$]
{\includegraphics[width=0.24\textwidth]{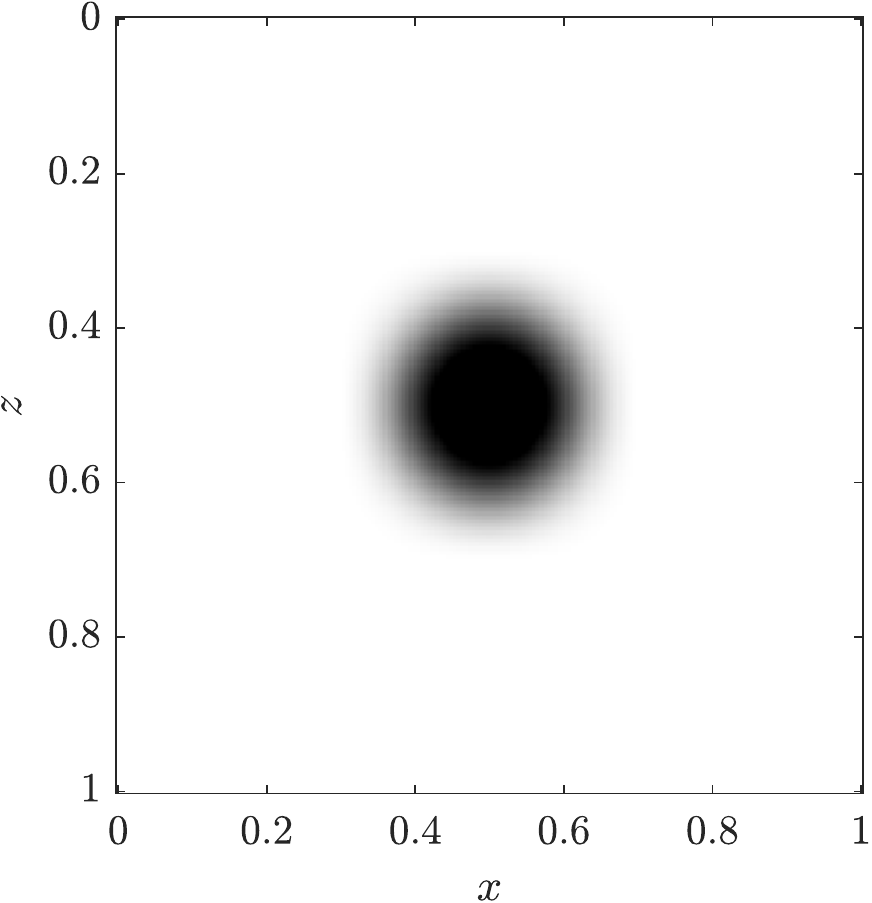}\label{fig:good true}}
\subfloat[$\fa = \fb = 0$]
{\includegraphics[width=0.24\textwidth]{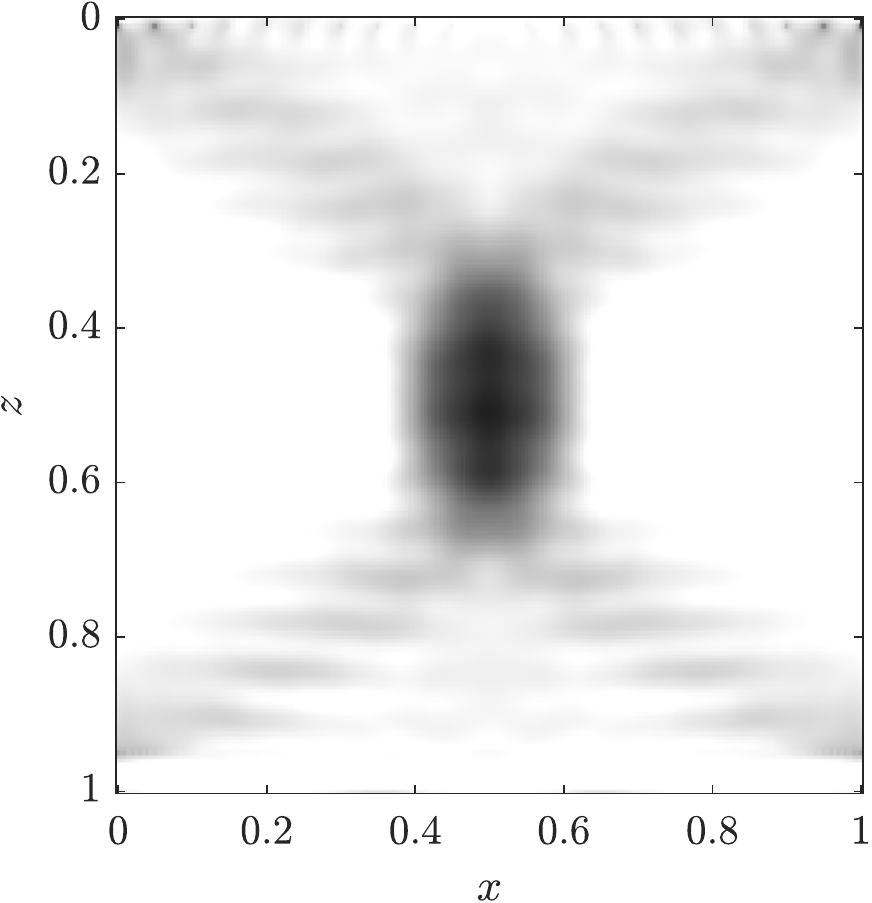}\label{fig:good a0}}
\subfloat[$\fa = 1$, $\fb = 0$]{\includegraphics[width=0.24\textwidth]{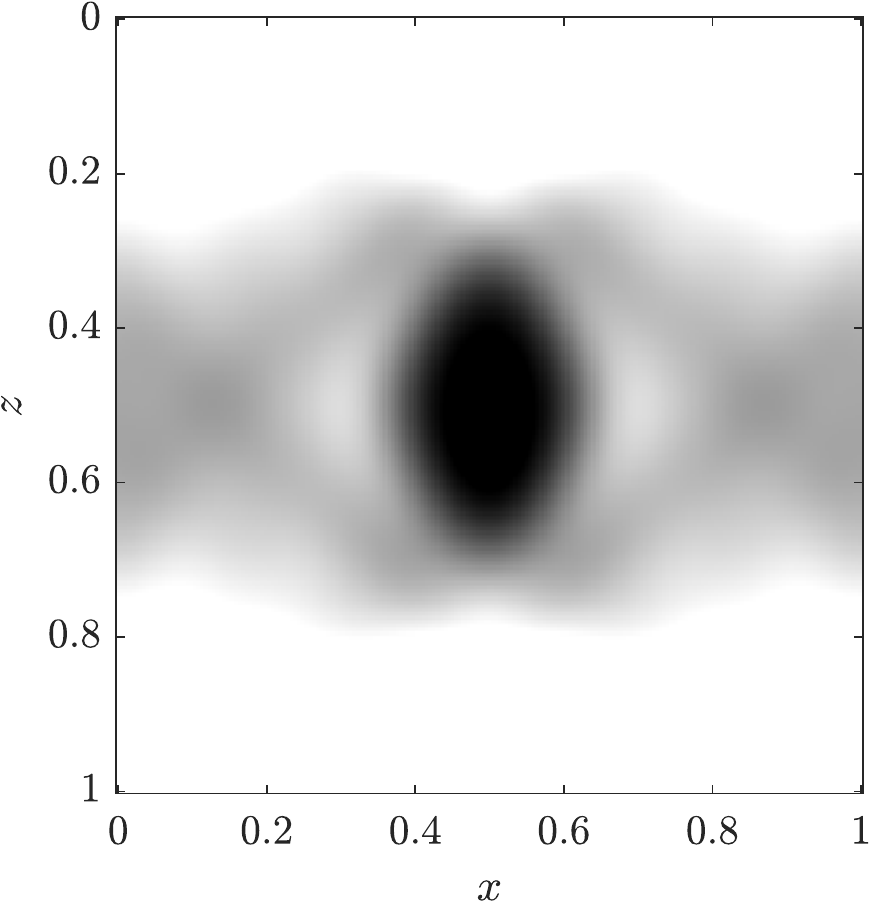}\label{fig:good a = 1}}
\subfloat[$\fa = 0$, $\fb = -1$]{\includegraphics[width=0.24\textwidth]{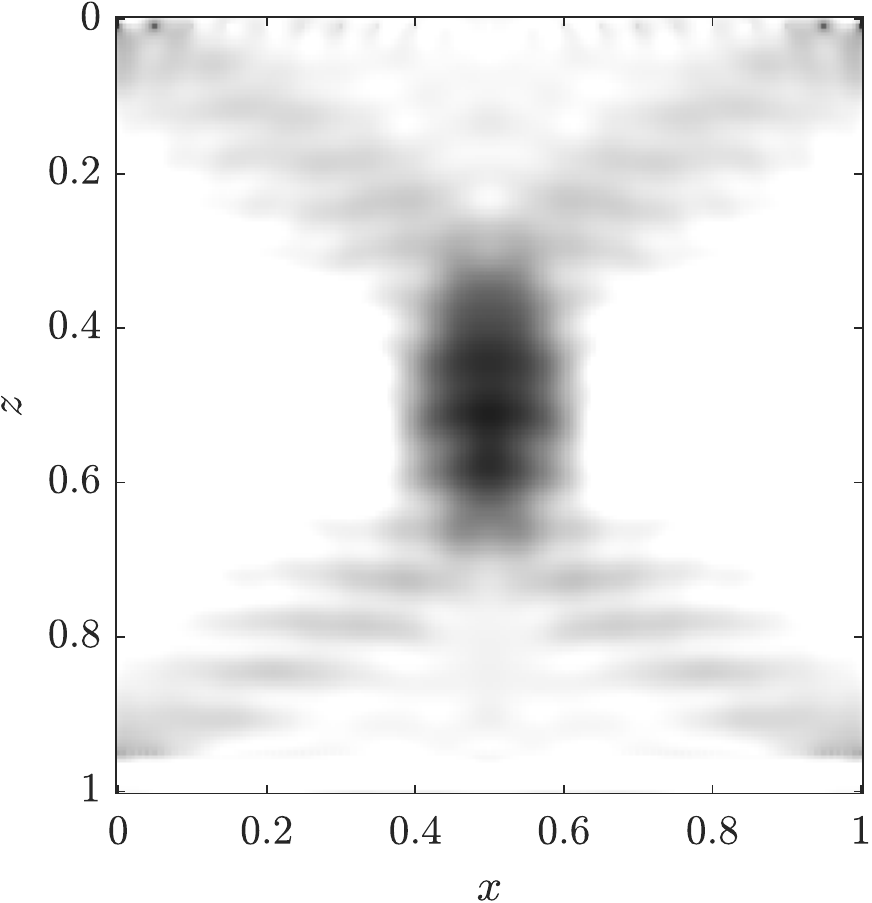}\label{fig:good b = -1}}
\caption{(Setting Two) (a): The true scatterer $\eta$. (b)-(d): The reconstructed results after $100$ iterations of L-BFGS algorithms under different $(\fa,\fb)$ choices.}
\label{fig:good recovery}
\end{figure}

\subsection{Inversion Examples}
Next, we present a few inversion examples based on the scalar wave equation~\eref{eq:helm-1}, where the forward problem is not too smoothing, and the phenomena of using different norms are easy to observe. 

We consider a square domain $\Omega = [0,1]^2$ discretized with the spatial spacing $dx = dz = 0.005$. The spatial unit is kilometer (km). We use point sources placed on the top of the domain at depth $z_s = 0.1$ and horizontal location $x_s \in \{0,0.05,0.1,\ldots, 1\}$. The measured data are recorded from the bottom of the domain at the depth $z_r = 0.95$ and horizontal location $x_r \in \{0,2dx,4dx,\ldots,1\}$. We use a single frequency $\omega = 21$~Hz and the background velocity $c_0 = 2.5$ km/s, so we have $k = \omega/c_0 = 8.4$; see~\eref{eq:helm-2}. 

We consider two different true scatterer $\eta$, which are shown in~\fref{fig:bad true} (Setting One) and~\fref{fig:good true} (Setting Two), respectively. Note that both scatterers are $C^\infty$, sharing a similar circular structure, but the one in~\fref{fig:good true} is a smoothed version of the scatterer in~\fref{fig:bad true} through a Gaussian filter. We choose these two examples to demonstrate the different radius of convergence for the inverse scattering series under different $(\fa,\fb)$ pairs. The inversion results under Setting One are shown in~\fref{fig:bad recovery} for three cases: $(\fa = 0, \fb = 0)$, $(\fa = 1, \fb = 0)$, and $(\fa = 0, \fb = -1)$. Similarly, inversion results under the same three cases for Setting Two are shown in the rest of~\fref{fig:good recovery}. 

From~\fref{fig:bad a0}, we observe that the less smooth scatterer shown in~\fref{fig:bad true} is \textit{outside} the radius of convergence by considering the inverse scattering as a map from $L^2$ to $L^2$. \Cref{fig:bad a0} is full of wrong features and shares few similarities with the ground truth.  In contrast, \fref{fig:bad true} is \textit{within} the radius of convergence under the setups of $(\fa = 0, \fb = -1$), and $(\fa = 1, \fb = 0)$ since~\fref{fig:bad a = 1} and~\fref{fig:bad b = -1} both recover the location and concentration of the scatterer, which are the main features of~\fref{fig:bad true}. We remark that this phenomenon matches our theoretical results in~\sref{section:inverse_scattering_series} and our follow-up discussions in~\sref{sec:discussions} where we conclude that both increasing $\fa$ and decreasing $\fb$ will increase the radius of convergence for the inverse scattering series.

\begin{figure}
\centering
\subfloat[Noisy data]
{\includegraphics[width=0.24\textwidth]{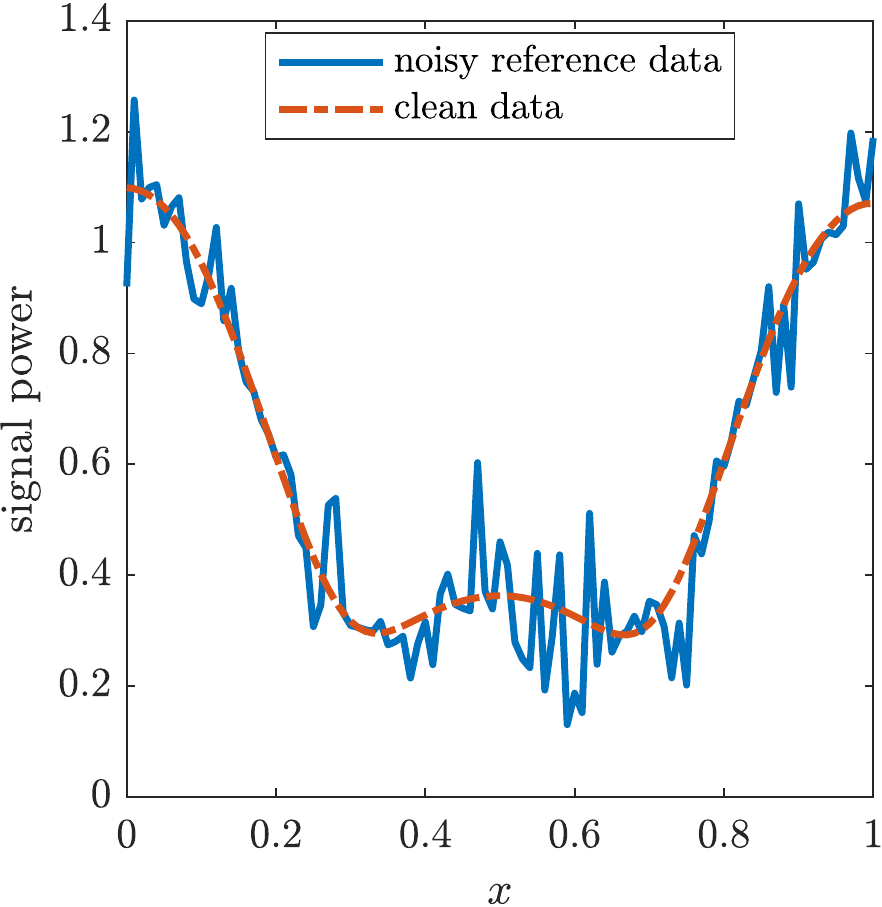}\label{fig:noisy data}}
\subfloat[$\fa = \fb = 0$]
{\includegraphics[width=0.24\textwidth]{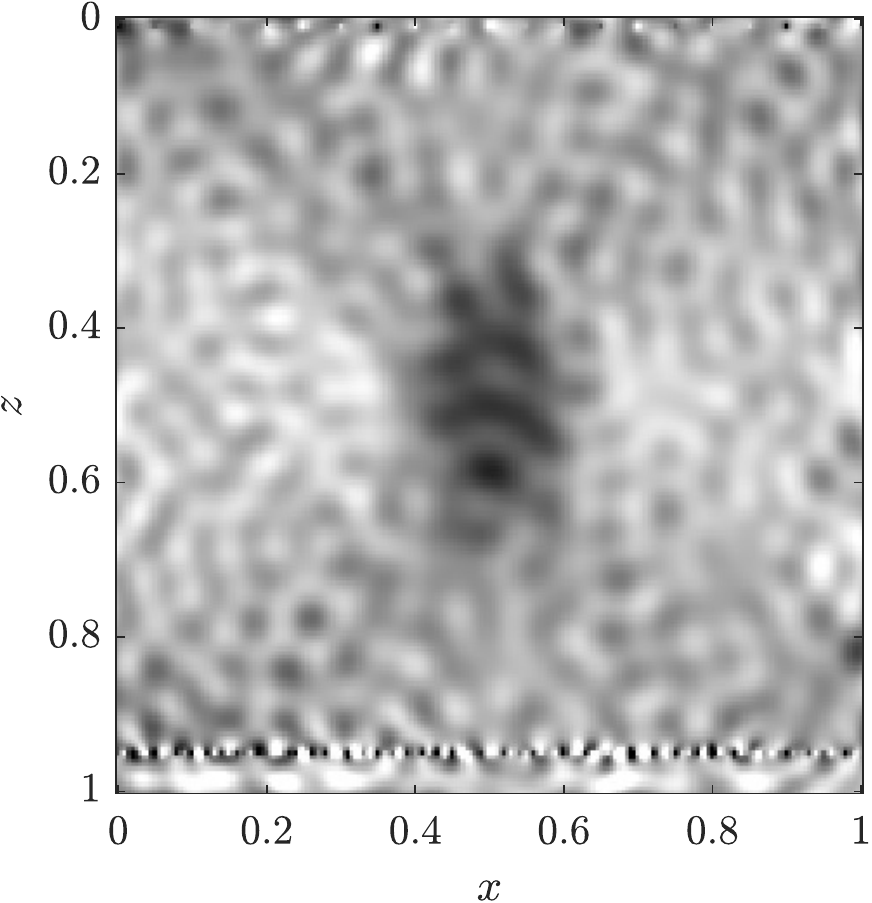}\label{fig:noisy a0}}
\subfloat[$\fa = 1$, $\fb = 0$]{\includegraphics[width=0.24\textwidth]{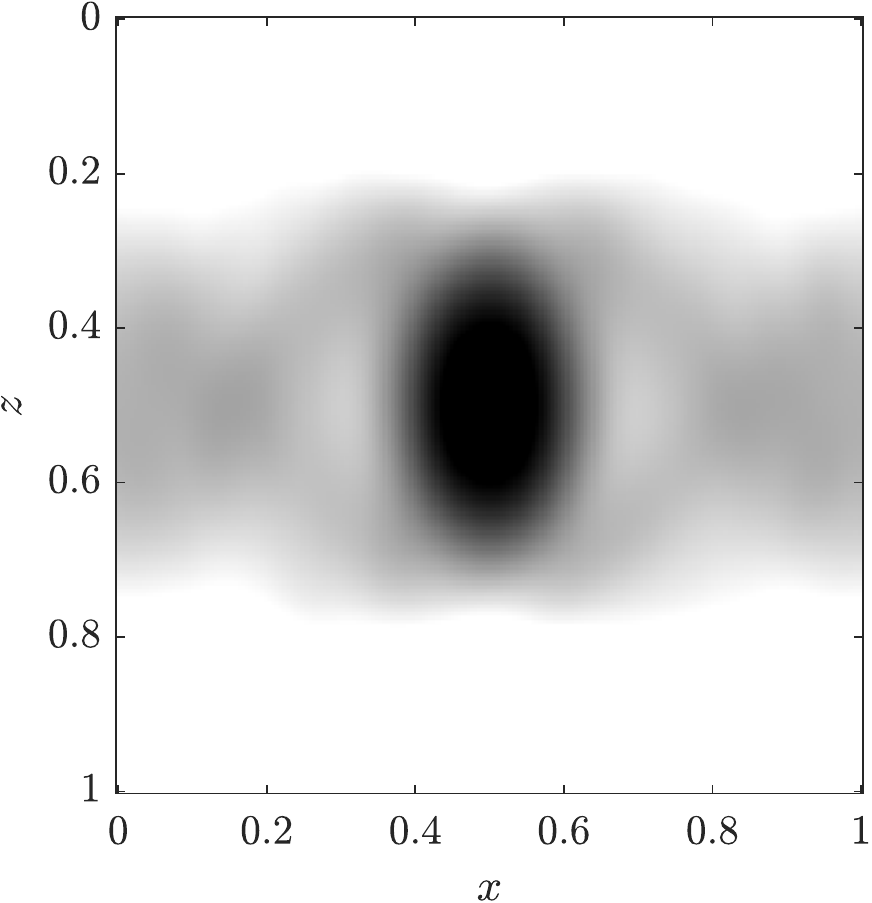}\label{fig:noisy a = 1}}
\subfloat[$\fa = 0$, $\fb = -1$]{\includegraphics[width=0.24\textwidth]{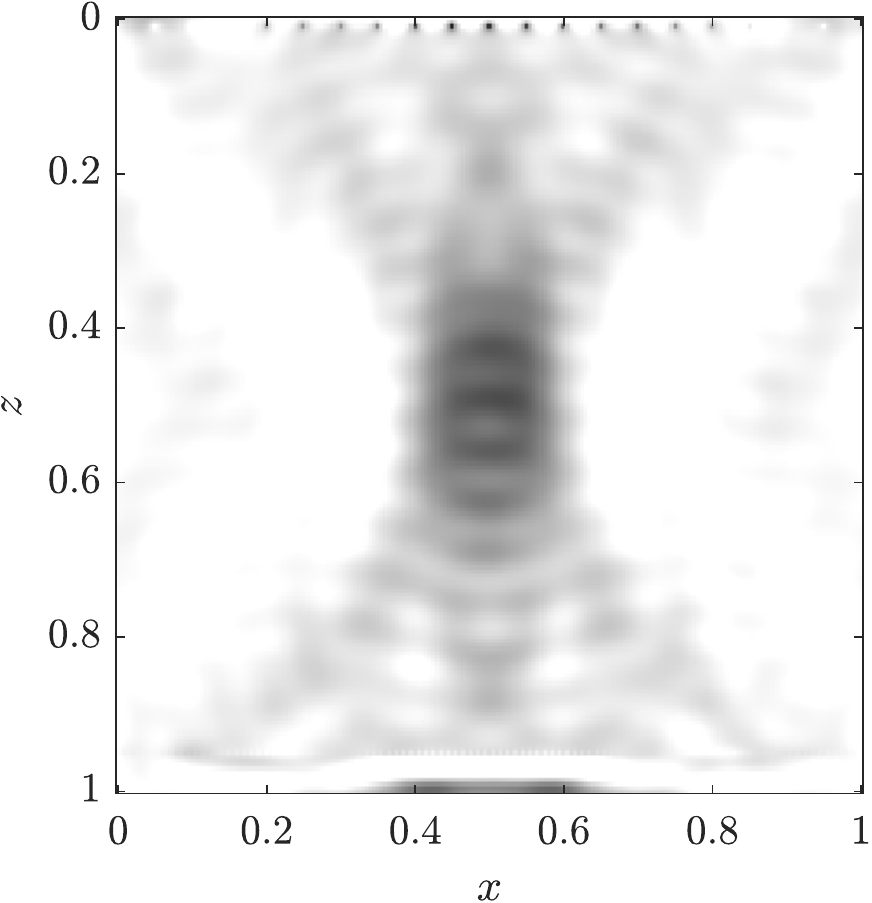}\label{fig:noisy b = -1}}
\caption{(Setting Two with noise) (a): Noisy data from one of the point sources. (b)-(d): The reconstructed results after $100$ iterations of L-BFGS algorithms under different $(\fa,\fb)$ choices while the reference data is polluted by white noise as shown in~(a).}
\label{fig:noisy recovery}
\end{figure}

For comparison, we also perform inversion where the ground truth is the smoother scatterer shown in~\fref{fig:good true}. When $\fa = \fb = 0$, the inverted result in~\fref{fig:good a0} also recovers the primary features of the ground truth, such as the location, concentration, and contrast of the scatterer. This shows that the $\eta$ shown in~\fref{fig:good true} is \textit{within} the radius of convergence for the case $\fa = \fb = 0$ while the one in~\fref{fig:bad true}  is not. The inversion result in~\fref{fig:good a = 1} is the closest to the ground truth after the same number of iterations since the choice $\fa = 1$ imposes an \textit{a priori} assumption on the regularity of the scatterer, we aim to invert, and meanwhile, our ground truth is indeed a smooth function. Note that we do not have an explicit regularization term here, but the choice of $\fb$ plays a role in enforcing an \textit{implicit} regularization for the inversion. On the other hand, the recovery in~\fref{fig:good b = -1} using $\fb = -1$ is visually more oscillatory than the case $\fb=0$ in~\fref{fig:good a0} while both set $\fa = 0$. This is because the squared $H^{-1}$ norm as the objective function has small weights on the high-wavenumber components of the data misfit and thus overlooked the differences, resulting in slow convergence on the high-wavenumber components of the parameter. We refer interested readers to~\cite{engquist2020quadratic,yang2021implicit,ren2022generalized} for a detailed discussion on the trade-offs of using the $H^s$ norms. We also note that in practice, varying $\fa$ and $\fb$ can improve the reconstructed results as it can ensure that $\eta$ is in the radius of convergence early on in optimization (large $\fa$ or small $\fb$) while later fall into regimes with better data sensitivity to improve resolution (small $\fa$ or large $\fb$).


Having illustrated the different radius of convergence, we show an example regarding the stability with respect to data noises. The setup remains the same as Setting Two except that the reference data is now polluted with white noise; see~\fref {fig:noisy data} for an illustration. The inversion results under the same three cases are shown in the rest of~\fref{fig:noisy recovery}. We do not have a regularization term in the objective function. The $L^2\mapsto L^2$ setting in~\fref{fig:bad a0} overfit the noise (see~\fref{fig:noisy a0}), while the a priori assumption that $\eta \in H^1$ prevented the reconstruction in the case $(\fa = 1, \fb = 0)$ from being polluted by the noisy data (see~\fref{fig:noisy a = 1}). Finally, the inversion result by embedding the noisy data into the weaker $H^{-1}$ space also helped mitigate the high-frequency noise due to the natural smoothing property of the weak norm (see~\fref{fig:noisy b = -1}).

\section{Conclusion}\label{sec:conclusion}
In this work, we generalize the results in~\cite{moskow2008convergence} and consider a different class of function spaces to analyze the inverse scattering series for Helmholtz and diffuse wave equations. Specifically, we analyze the convergence property of the forward scattering series that maps the parameter in $H^\fa$ to the scattering data in $H^\fb$ and the convergence, stability, and approximation error of the inverse scattering series. Since we express all the results and constants explicitly in terms of $\fa$ and $\fb$, we could analyze how the choices of $H^\fa$ and $H^\fb$ for the parameter and the scattering data, respectively, impact certain properties of the inverse scattering problem. For example, we can increase the radius of convergence for the inverse scattering series by using a stronger function space $H^\fa$ for the parameter while using a weaker metric space $H^\fb$ for the scattering data, which corresponds to better convexity when using weaker $H^s$ norm as the objective function in practice~\cite{yang2021implicit}. Also, we can obtain a smaller stability constant by considering a weaker space for the parameter $H^\fa$, but it is at the cost of obtaining a reconstruction with less regularity.

\section*{Acknowledgements}
This material is based upon work supported by the National Science Foundation under Award Number DMS-1913129. S.~Mahankali acknowledges the PRIMES-USA program for making the collaboration possible. Y.~Yang acknowledges support from Dr.~Max R\"ossler, the Walter Haefner Foundation, and the ETH Z\"urich Foundation. We would like to thank the anonymous referees for their useful comments that helped us improve the quality of this work.

\appendix

\section{Proofs of Lemmas and Theorems in Section 4}\label{sec:appendix}
Here, we provide the proof of \Cref{lem:new_inverse_operator_bound}  and \Cref{thm:inverse_approx_error}. The proofs are based on results and proof techniques in~\cite[Sec.~3]{moskow2008convergence}. We include the proofs simply for completeness, as we explicitly state the values of all constants in~\Sref{section:inverse_scattering_series}.
 
\subsection{Proof of \Cref{lem:new_inverse_operator_bound}}

We provide a bound on the quantity $\|\mathcal{K}_j\|_{\fba}$ in \Cref{lem:new_inverse_operator_bound}. 

\begin{proof}[Proof of \Cref{lem:new_inverse_operator_bound}]
Consider $j \ge 2$. From \eref{eq:curly_K_j_definition}, we first have the bound
\begin{eqnarray} 
    \|\mathcal{K}_j\|_{\fba} &\le \left(\sum_{m=1}^{j-1}\|\mathcal{K}_m\|_{\fba}\sum_{i_1 + \cdots + i_m = j}\|K_{i_1}\|_{\ha \to \hb}\cdots\|K_{i_m}\|_{\ha \to \hb}\right) \|\mathcal{K}_1\|_{\fba}^j \nonumber \\
    &\le \|\mathcal{K}_1\|_{\fba}^j\left(\sum_{m=1}^{j-1}\|\mathcal{K}_m\|_{\fba}\sum_{i_1 + \cdots + i_m = j}\nu_{\fab}\mu_{\fab}^{i_1-1}\cdots\nu_{\fab}\mu_{\fab}^{i_m-1}\right), \label{eq:curly_K_j_first_bound} 
\end{eqnarray}
where we have used \Cref{lem:s_1_to_s_2_operator_bound} to obtain the second inequality. Using~\eref{eq:curly_K_j_first_bound}, the rest of the proof follows similarly as in~\cite[Lemma~2.1]{moskow2019inverse}.
\end{proof}

\subsection{Proof of \Cref{thm:inverse_approx_error}}

While the inverse scattering series converges under the conditions of \Cref{thm:inverse_series_convergence}, its limit is not equal to $\eta$ in general. In \Cref{thm:inverse_approx_error}, we give a bound on the distance between this limit and $\eta$.

\begin{proof}[Proof of \Cref{thm:inverse_approx_error}]
The first part of the proof uses the following derivations from~\cite[Theorem~2.1]{moskow2019inverse}. Since
$\|\mathcal{K}_1\|_{\fba} < 1/(\mu_{\fab} + \nu_{\fab})$ and $\|\mathcal{K}_1\phi\|_{\ha\times L^2} < 1/(\mu_{\fab} + \nu_{\fab})$, the inverse scattering series
\begin{equation}\label{eq:inverse_series_formula}
    \widetilde{\eta} = \sum_{j}\mathcal{K}_j\phi\otimes\cdots\otimes\phi
\end{equation}
converges. Also, since $\mathcal{M} < 1/(\mu_{\fab} + \nu_{\fab})$, the forward scattering series
\[
    \phi = \sum_{j}K_j\eta\otimes\cdots\otimes\eta
\]
converges as well. Substituting the forward scattering series for $\phi$ into \eref{eq:inverse_series_formula} gives
\[
    \widetilde{\eta} = \sum_{j}\widetilde{\mathcal{K}}_j\eta\otimes\cdots\otimes\eta
\]
where $\widetilde{\mathcal{K}}_1 = \mathcal{K}_1K_1$ and for $j \ge 2$,
\[
    \widetilde{\mathcal{K}}_j = \left(\sum_{m=1}^{j-1}\mathcal{K}_m\sum_{i_1 + \cdots + i_m = j}K_{i_1}\otimes\cdots\otimes K_{i_m}\right) + \mathcal{K}_jK_1\otimes\cdots\otimes K_1,
\]
as defined in~\cite[Eqn.~(67)]{moskow2008convergence}.
Using the formula for $\mathcal{K}_j$ in \eref{eq:curly_K_j_definition}, we get 
\[
    \widetilde{\mathcal{K}}_j = \left(\sum_{m=1}^{j-1}\mathcal{K}_m\sum_{i_1 + \cdots + i_m = j}K_{i_1}\otimes\cdots\otimes K_{i_m}\right)(I - \mathcal{K}_1K_1\otimes\cdots\otimes\mathcal{K}_1K_1).
\]
Since $\eta - \widetilde{\eta} = \eta - \mathcal{K}_1K_1\eta - \widetilde{\mathcal{K}}_2\eta\otimes\eta + \cdots$, we get
\[
    \|\eta - \widetilde{\eta}\|_{\ha}
    \le \|\eta - \mathcal{K}_1K_1\eta\|_{\ha} + \sum_{j=2}^{\infty}\|\widetilde{\mathcal{K}}_j\eta\otimes\cdots\otimes\eta\|_{\ha}
\]
using the triangle inequality. Next, consider ${\|\widetilde{\mathcal{K}}_j\eta\otimes\cdots\otimes\eta\|_{\ha}}$, for $j \ge 2$. We first get
\begin{eqnarray*}
    \|\widetilde{\mathcal{K}}_j\eta\otimes\cdots\otimes\eta\|_{\ha} &\le \sum_{m=1}^{j-1} \|\mathcal{K}_m\|_{\fba} \sum_{i_1 + \cdots + i_m = j} \|K_{i_1}\|_{\ha \to \hb}\cdots\|K_{i_m}\|_{\ha \to \hb}&\\
    &\cdot \|\eta\otimes\cdots\otimes\eta - \mathcal{K}_1K_1\eta\otimes\cdots\otimes\mathcal{K}_1K_1\eta\|_{\ha}.&
\end{eqnarray*}
Letting $\psi = \eta - \mathcal{K}_1K_1\eta$ and using~\cite[Eqn.~(29)]{moskow2019inverse} gives
\[
    \|\eta\otimes\cdots\otimes\eta - \mathcal{K}_1K_1\eta\otimes\cdots\otimes\mathcal{K}_1K_1\eta\|_{\ha} \le j\mathcal{M}^{j-1}\|\psi\|_{\ha}.
\]
Our proof is distinct from the literature after this point. From the proof of \Cref{lem:new_inverse_operator_bound}, 
\[ 
    \sum_{m=1}^{j-1} \|\mathcal{K}_m\|_{\fba} \sum_{i_1 + \cdots + i_m = j} \|K_{i_1}\|_{\ha \to \hb}\cdots\|K_{i_m}\|_{\ha \to \hb} \le C(\mu_{\fab} + \nu_{\fab})^j
\]
since the left-hand side of this equation is the same as the right-hand side of \eref{eq:curly_K_j_first_bound}, except with the factor of $\|\mathcal{K}_1\|_{\fba}^j$ removed.
This means, for $j \ge 2$, we have
\[
    \|\widetilde{\mathcal{K}}_j\eta\otimes\cdots\otimes\eta\|_{\ha} \le C\|\psi\|_{\ha}j\mathcal{M}^{j-1}(\mu_{\fab} + \nu_{\fab})^j,
\]
which gives
\begin{eqnarray}
    \|\eta - \widetilde{\eta}\|_{\ha}
    &\le \|\psi\|_{\ha} + C\|\psi\|_{\ha}\sum_{j=2}^{\infty}j\mathcal{M}^{j-1}(\mu_{\fab} + \nu_{\fab})^j \nonumber \\
    &\le C^*\|\psi\|_{\ha}\sum_{j=1}^{\infty}j\mathcal{M}^{j-1}(\mu_{\fab} + \nu_{\fab})^j, \label{eq:sharper_ineq}
\end{eqnarray}
and this series converges since $(\mu_{\fab} + \nu_{\fab})\mathcal{M} < 1$. Computing this sum gives
\begin{equation} \label{eq:sharper_ineq2}
    \|\eta - \widetilde{\eta}\|_{\ha} \le
    \frac{C^*(\mu_{\fab} + \nu_{\fab})\|\eta - \mathcal{K}_1K_1\eta\|_{\ha}}{(1 - (\mu_{\fab} + \nu_{\fab})\mathcal{M})^2} = C_{\fab}\|\eta - \mathcal{K}_1K_1\eta\|_{\ha}.
\end{equation}
We can split the approximation error into two terms with the triangle inequality:
\begin{equation*} \fl
    \left|\left|\eta - \sum_{j=1}^N \mathcal{K}_j\phi\otimes\cdots\otimes\phi\right|\right|_{\ha}  \le \left|\left|\eta - \sum_{j} \widetilde{\mathcal{K}}_j\eta\otimes\cdots\otimes\eta\right|\right|_{\ha}+  \left|\left|\widetilde{\eta} - \sum_{j=1}^N \mathcal{K}_j\phi\otimes\cdots\otimes\phi\right|\right|_{\ha}.
\end{equation*}
Using \Cref{thm:inverse_series_convergence}, we finally obtain
\begin{equation*} \fl
    \left|\left|\eta - \sum_{j=1}^N \mathcal{K}_j\phi\otimes\cdots\otimes\phi\right|\right|_{\ha}  \le C_{\fab}\| \eta - \mathcal{K}_1K_1\eta\|_{\ha} + C\frac{\left[(\mu_{\fab} + \nu_{\fab})\|\mathcal{K}_1\|_{\fba}\|\phi\|_{\hb}\right]^{N+1}}{1 - (\mu_{\fab} + \nu_{\fab})\|\mathcal{K}_1\|_{\fba}\|\phi\|_{\hb}}
\end{equation*}
where $C$ is defined in \eref{eq:C_original_definition}.
\end{proof}

\section*{References}
\bibliography{references}

\begin{thebibliography}{10}

\bibitem{arbogast1999methods}
Todd Arbogast and Jerry~L Bona.
\newblock Methods of applied mathematics.
\newblock {\em Department of Mathematics, University of Texas}, 2008, 1999.

\bibitem{JS4}
Simon Arridge, Shari Moskow, and John~C Schotland.
\newblock Inverse {Born} series for the calderon problem.
\newblock {\em Inverse Problems}, 28(3):035003, 2012.

\bibitem{Bao2015}
Gang Bao, Peijun Li, Junshan Lin, and Faouzi Triki.
\newblock {Inverse scattering problems with multi-frequencies}.
\newblock {\em Inverse Problems}, 31(9):093001, {S}ep 2015.

\bibitem{barcelo2021fourier}
JA~Barcel{\'o}, M~Folch-Gabayet, T~Luque, S~P{\'e}rez-Esteva, and MC~Vilela.
\newblock {The Fourier extension operator of distributions in Sobolev spaces of
  the sphere and the Helmholtz equation}.
\newblock {\em Proceedings of the Royal Society of Edinburgh Section A:
  Mathematics}, 151(6):1768--1789, 2021.

\bibitem{barcelo2020characterization}
Juan~Antonio Barcel{\'o}, Teresa Luque, and Salvador P{\'e}rez-Esteva.
\newblock {Characterization of Sobolev spaces on the sphere}.
\newblock {\em Journal of Mathematical Analysis and Applications},
  491(1):124240, 2020.

\bibitem{bardsley2014restarted}
Patrick Bardsley and Fernando~Guevara Vasquez.
\newblock Restarted inverse {Born} series for the {S}chr{\"o}dinger problem
  with discrete internal measurements.
\newblock {\em Inverse Problems}, 30(4):045014, 2014.

\bibitem{Barnett2017}
Alex Barnett, Leslie Greengard, Andras Pataki, and Marina Spivak.
\newblock {Rapid Solution of the Cryo-EM Reconstruction Problem by Frequency
  Marching}.
\newblock {\em SIAM Journal on Imaging Sciences}, 10(3):1170--1195, Jan 2017.

\bibitem{Brossier2010}
Romain Brossier, St{\'e}phane Operto, and Jean Virieux.
\newblock Which data residual norm for robust elastic frequency-domain full
  waveform inversion?
\newblock {\em Geophysics}, 75(3):R37--R46, 2010.

\bibitem{bunks1995multiscale}
Carey Bunks, Fatimetou~M Saleck, S~Zaleski, and G~Chavent.
\newblock Multiscale seismic waveform inversion.
\newblock {\em Geophysics}, 60(5):1457--1473, 1995.

\bibitem{calderMY10}
Jeff Calder, A~Mansouri, and Anthony Yezzi.
\newblock Image sharpening via {S}obolev gradient flows.
\newblock {\em SIAM Journal on Imaging Sciences}, 3(4):981--1014, 2010.

\bibitem{JS7}
Francis~J Chung, Anna~C Gilbert, Jeremy~G Hoskins, and John~C Schotland.
\newblock Optical tomography on graphs.
\newblock {\em Inverse Problems}, 33(5):055016, 2017.

\bibitem{colton1998inverse}
David~L Colton and Rainer Kress.
\newblock {\em Inverse acoustic and electromagnetic scattering theory},
  volume~93.
\newblock Springer, 1998.

\bibitem{dunlop2020stability}
Matthew~M Dunlop and Yunan Yang.
\newblock {Stability of Gibbs posteriors from the Wasserstein loss for Bayesian
  full waveform inversion}.
\newblock {\em SIAM/ASA Journal on Uncertainty Quantification},
  9(4):1499--1526, 2021.

\bibitem{engquist2020quadratic}
Bj{\"o}rn Engquist, Kui Ren, and Yunan Yang.
\newblock The quadratic {W}asserstein metric for inverse data matching.
\newblock {\em Inverse Problems}, 36(5):055001, 2020.

\bibitem{ren2022generalized}
Bj{\"o}rn Engquist, Kui Ren, and Yunan Yang.
\newblock A generalized weighted optimization method for computational learning
  and inversion.
\newblock In {\em International Conference on Learning Representations}, 2022.

\bibitem{federer2014geometric}
Herbert Federer.
\newblock {\em Geometric measure theory}.
\newblock Springer, 2014.

\bibitem{Fichtner2013}
Andreas Fichtner, Jeannot Trampert, Paul Cupillard, Erdinc Saygin, Tuncay
  Taymaz, Yann Capdeville, and Antonio Villase{\~{n}}or.
\newblock {Multiscale full waveform inversion}.
\newblock {\em Geophysical Journal International}, 194(1):534--556, Jul 2013.

\bibitem{hoskins2022analysis}
Jeremy~G Hoskins and John~C Schotland.
\newblock {Analysis of the inverse {Born} series: an approach through geometric
  function theory}.
\newblock {\em Inverse Problems}, 38(7):074001, 2022.

\bibitem{jacobs2019solving}
Matt Jacobs, Flavien L{\'e}ger, Wuchen Li, and Stanley Osher.
\newblock Solving large-scale optimization problems with a convergence rate
  independent of grid size.
\newblock {\em SIAM Journal on Numerical Analysis}, 57(3):1100--1123, 2019.

\bibitem{jeong2017comparison}
Gangwon Jeong, Jongha Hwang, and Dong-Joo Min.
\newblock Comparison of weighting techniques for acoustic full waveform
  inversion.
\newblock {\em Journal of Applied Geophysics}, 147:16--27, 2017.

\bibitem{Kamei2014}
Rie Kamei, R~Gerhard Pratt, and Takeshi Tsuji.
\newblock {Misfit functionals in Laplace--Fourier domain waveform inversion,
  with application to wide-angle ocean bottom seismograph data}.
\newblock {\em Geophysical Prospecting}, 62:1054--1074, 2014.

\bibitem{JS5}
Kimberly Kilgore, Shari Moskow, and John~C Schotland.
\newblock Inverse {Born} series for scalar waves.
\newblock {\em Journal of Computational Mathematics}, pages 601--614, 2012.

\bibitem{JS3}
Kimberly Kilgore, Shari Moskow, and John~C Schotland.
\newblock Convergence of the {{Born}} and inverse {Born} series for
  electromagnetic scattering.
\newblock {\em Applicable analysis}, 96(10):1737--1748, 2017.

\bibitem{leoni2017first}
Giovanni Leoni.
\newblock {\em A first course in Sobolev spaces}.
\newblock American Mathematical Soc., 2017.

\bibitem{JS6}
Manabu Machida and John~C Schotland.
\newblock Inverse {Born} series for the radiative transport equation.
\newblock {\em Inverse Problems}, 31(9):095009, 2015.

\bibitem{mandelis2000diffusion}
Andreas Mandelis.
\newblock Diffusion waves and their uses.
\newblock {\em Physics today}, 53(8):29--34, 2000.

\bibitem{mandelis2013diffusion}
Andreas Mandelis.
\newblock {\em Diffusion-wave fields: mathematical methods and Green
  functions}.
\newblock Springer Science \& Business Media, 2013.

\bibitem{JS1}
Vadim~A Markel, Joseph~A O’Sullivan, and John~C Schotland.
\newblock Inverse problem in optical diffusion tomography. iv. nonlinear
  inversion formulas.
\newblock {\em JOSA A}, 20(5):903--912, 2003.

\bibitem{moskow2008convergence}
Shari Moskow and John~C Schotland.
\newblock Convergence and stability of the inverse scattering series for
  diffuse waves.
\newblock {\em Inverse Problems}, 24(6):065005, 2008.

\bibitem{JS2}
Shari Moskow and John~C Schotland.
\newblock {Numerical studies of the inverse {Born} series for diffuse waves}.
\newblock {\em Inverse Problems}, 25(9):095007, 2009.

\bibitem{moskow2019inverse}
Shari Moskow and John~C Schotland.
\newblock Inverse {Born} series.
\newblock {\em The Radon Transform, The First 100 Years and Beyond}, 22, 2019.

\bibitem{necas2011direct}
Jindrich Necas.
\newblock {\em Direct methods in the theory of elliptic equations}.
\newblock Springer Science \& Business Media, 2011.

\bibitem{neuberger2009sobolev}
John Neuberger.
\newblock {\em Sobolev gradients and differential equations}.
\newblock Springer Science \& Business Media, 2009.

\bibitem{nurbekyan2022efficient}
Levon Nurbekyan, Wanzhou Lei, and Yunan Yang.
\newblock Efficient natural gradient descent methods for large-scale
  optimization problems.
\newblock {\em arXiv preprint arXiv:2202.06236}, 2022.

\bibitem{Oh2013}
Ju-Won Oh and Dong-Joo Min.
\newblock {Weighting technique using backpropagated wavefields incited by
  deconvolved residuals for frequency-domain elastic full waveform inversion}.
\newblock {\em Geophysical Journal International}, 194:322--347, 2013.

\bibitem{otto2000generalization}
Felix Otto and C{\'e}dric Villani.
\newblock Generalization of an inequality by {T}alagrand and links with the
  logarithmic {S}obolev inequality.
\newblock {\em Journal of Functional Analysis}, 173(2):361--400, 2000.

\bibitem{papadakis2014optimal}
Nicolas Papadakis, Gabriel Peyr{\'e}, and Edouard Oudet.
\newblock Optimal transport with proximal splitting.
\newblock {\em SIAM Journal on Imaging Sciences}, 7(1):212--238, 2014.

\bibitem{peyre2018comparison}
R{\'e}mi Peyre.
\newblock Comparison between ${W}_2$ distance and $\dot{H}^{-1}$ norm, and
  localization of {W}asserstein distance.
\newblock {\em ESAIM: Control, Optimisation and Calculus of Variations},
  24(4):1489--1501, 2018.

\bibitem{rektorys2012variational}
Karel Rektorys.
\newblock {\em Variational methods in mathematics, science and engineering}.
\newblock Springer Science \& Business Media, 2012.

\bibitem{ryan2002introduction}
Raymond~A Ryan.
\newblock {\em {Introduction to tensor products of Banach spaces}}, volume~73.
\newblock Springer, 2002.

\bibitem{sundaramoorthi2007sobolev}
Ganesh Sundaramoorthi, Anthony Yezzi, and Andrea~C Mennucci.
\newblock Sobolev active contours.
\newblock {\em International Journal of Computer Vision}, 73(3):345--366, 2007.

\bibitem{wieczorek2018shtools}
Mark~A Wieczorek and Matthias Meschede.
\newblock Shtools: Tools for working with spherical harmonics.
\newblock {\em Geochemistry, Geophysics, Geosystems}, 19(8):2574--2592, 2018.

\bibitem{yang2018application}
Yunan Yang, Bj{\"o}rn Engquist, Junzhe Sun, and Brittany~F Hamfeldt.
\newblock Application of optimal transport and the quadratic {W}asserstein
  metric to full-waveform inversion.
\newblock {\em Geophysics}, 83(1):R43--R62, 2018.

\bibitem{yang2021implicit}
Bowen Zhu, Jingwei Hu, Yifei Lou, and Yunan Yang.
\newblock {Implicit regularization effects of the Sobolev norms in image
  processing}.
\newblock {\em arXiv preprint arXiv:2109.06255}, 2021.

\end{thebibliography}
\bibliographystyle{plain}
\end{document}